# On quasi-reductive group schemes

Gopal Prasad[1]
University of Michigan
and
Jiu-Kang Yu[2]
Purdue University

## 1 Introduction

In [SGA3], the following remarkable theorem about tori is proved:

**1.1 Theorem** [SGA3, Exp. X, Théorème 8.8] *Let $\mathcal{T}$ be a commutative flat group scheme, separated of finite type over a noetherian scheme S, with connected affine fibers. Let $s \in S$, $\bar{s}$ a geometric point over s, and suppose*

- *the reduced subscheme $(\mathcal{T}_{\bar{s}})_{\mathrm{red}}$ of the geometric fiber $\mathcal{T}_{\bar{s}}$ is a torus; and*
- *there exists a generization t of s (i.e. the closure of $\{t\}$ contains s) such that $\mathcal{T}_t$ is smooth over $\kappa(t)$, the residue field of t.*

*Then there exists an open neighborhood U of s such that $\mathcal{T}|U$ is a torus over U.*

The aim of this paper is to prove an analogous result where "torus" is replaced by "reductive group". As usual, the key point is to treat the case where the base scheme $S$ is the spectrum of a discrete valuation ring (henceforth to be called a DVR for brevity) $R$. In this case, our result answers a question of Kari Vilonen. It turns out that a direct analogue of Theorem 1.1 fails to hold in some cases if the residue field of $R$ is of characteristic 2, but the group schemes for which the theorem fails can be classified over a strictly henselian $R$.

To state our results, let $R$ be a DVR, $\pi$ a uniformizer of $R$ and $\kappa$ the residue field. Let $K = \mathrm{Frac}\, R$. We will call a group scheme $\mathcal{G}$ over $R$ *quasi-reductive* (this is unrelated to the notion of *quasi-reductive algebraic groups* over a non-perfect field introduced in [BT2, 1.1.12]) if

1) $\mathcal{G}$ is affine and flat over $R$,
2) the generic fiber $\mathcal{G}_K := \mathcal{G} \otimes_R K$ is connected and smooth over $K$,
3) the special fiber $\mathcal{G}_\kappa := \mathcal{G} \otimes_R \kappa$ is of finite type over $\kappa$ and the neutral component $(\mathcal{G}_{\bar\kappa})^\circ_{\mathrm{red}}$ of the reduced geometric special fiber $(\mathcal{G}_{\bar\kappa})_{\mathrm{red}}$ is a reductive group of dimension $= \dim \mathcal{G}_K$.

**1.2 Theorem** *Let $\mathcal{G}$ be a quasi-reductive group scheme over R. Then*

(a) *$\mathcal{G}$ is of finite type over R;*
(b) *the generic fiber $G := \mathcal{G}_K$ is reductive;*

---

[1]partially supported by an NSF grant
[2]partially supported by an NSF grant, a Sloan fellowship, and the IHES

(c) *the special fiber $\mathcal{G}_\kappa$ is connected;*

*In addition, $\mathcal{G}$ is a reductive group over $R$ if at least one of the following holds:*

(i) $\operatorname{char}\kappa \neq 2$;
(ii) *the type of $G \otimes \bar{K}$ is the same as that of $(\mathcal{G}_{\bar\kappa})^\circ_{\mathrm{red}}$;*
(iii) *no normal algebraic subgroup of $G \otimes \bar{K}$ is isomorphic to $\mathrm{SO}_{2n+1}$ for $n \geqslant 1$.*

We recall [SGA3, Exp. XIX, 2.7] that a *reductive group over $S$* is a smooth affine group scheme over $S$ such that all the geometric fibers are connected reductive algebraic groups. By the *type* of a connected split reductive algebraic group over a field, we mean the isomorphism class of its associated root datum $(X, \Phi, X^\vee, \Phi^\vee)$ [SGA3, Exp. XXII, 2.6]. If $R$ is strictly henselian, the isomorphism class of a reductive group over $R$ is determined by its type (i.e. the type of either $\mathcal{G}_K$ or $\mathcal{G}_\kappa$, which are the same [SGA3, Exp. XXII, Proposition 2.8, Exp. XXIII, Corollaire 5.2]).

We remark that over the base $\operatorname{Spec} R$, even in the case of tori, Theorem 1.2 is slightly stronger than Theorem 1.1 in that we do not assume that $\mathcal{G}$ is of finite type over $R$ (and we only impose conditions on $(\mathcal{G}_{\bar\kappa})^\circ_{\mathrm{red}}$). This generality is also required by Vilonen's question. Notice that if $\mathcal{G}$ is of finite type over $R$, then $\mathcal{G}_\kappa$ has the same dimension as $G$ by [EGA, IV, Lemme 14.3.10]. In §7, we will provide examples which show that without the condition $\dim \mathcal{G}_\kappa = \dim \mathcal{G}_K$ imposed in the definition of quasi-reductive group schemes, the preceding theorem is false.

In addition, we have the following:

**1.3 Corollary** *Let $\mathcal{G}$ be a reductive group over $R$ and $G$ be its generic fiber. Assume either*

- $\operatorname{char}\kappa \neq 2$; *or*
- *no normal algebraic subgroup of $G \otimes \bar{K}$ is isomorphic to $\mathrm{SO}_{2n+1}$ for $n \geqslant 1$.*

*Let $\phi : \mathcal{G} \to \mathcal{H}$ be a morphism of affine group schemes of finite type over $R$ such that $\phi_K : \mathcal{G}_K \to \mathcal{H}_K$ is a closed immersion. Then $\phi$ is a closed immersion.*

We remark that [SGA3, Exp. XVI, Proposition 1.5] implies the following statement: *Let $\phi : \mathcal{G} \to \mathcal{H}$ be a morphism of affine group schemes of finite type over $R$ such that $\mathcal{G}$ is reductive, and both $\phi_K : \mathcal{G}_K \to \mathcal{H}_K$ and $\phi_\kappa : \mathcal{G}_\kappa \to \mathcal{H}_\kappa$ are closed immersions, then $\phi$ is a closed immersion.* The above corollary shows that the hypothesis on $\phi_\kappa$ can be eliminated provided a restriction on $K$ or $G$ is imposed (see 7.1; we note here that an incorrect version of Corollary 1.3 was given in [V, Proposition 3.1.2.1c]). This result, together with the existence of Bruhat–Tits schemes corresponding to parahoric subgroups, has the following remarkable consequence:

**1.4 Corollary** *Assume that $R$ is strictly henselian and $\kappa$ is algebraically closed. Let $\phi : G \hookrightarrow G'$ be an inclusion of connected reductive algebraic groups over $K$, $P$ a hyperspecial parahoric subgroup of $G(K)$, $P'$ an arbitrary parahoric subgroup of $G'(K)$ such that $\phi(P) \subset P'$. Assume that either the characteristic of $\kappa$ is not 2, or no normal algebraic subgroup of $G \otimes \bar{K}$ is isomorphic to $\mathrm{SO}_{2n+1}$ for $n \geqslant 1$. Then any function $f \in K[G]$ such that $f(P) \subset R$ is the restriction of a function $f' \in K[G']$ such that $f'(P') \subset R$.*

An analogue of Theorem 1.2 over a general noetherian base scheme is given in §6. There, we will also prove the following:

**1.5 Proposition** *Let S be a Dedekind scheme and $\mathcal{G}$ an affine flat group scheme over S such that all the fibers are reductive groups of the same dimension. Then $\mathcal{G}$ is a reductive group over S.*

These results have been used by Mirković and Vilonen to give a geometric interpretation of the dual group [MV]. We note here that their group scheme arises from geometry via a Tannakian formalism and they do not know a priori if it is of finite type; however, the results of this paper show that their group scheme is of finite type.

We will now summarize the content of this paper. Sections 2–5 are devoted to the proof of Theorem 1.2 and Corollary 1.3. In Section 6, using Theorem 1.2, we prove results about group schemes over Dedekind schemes and group schemes of finite type over more general noetherian schemes. In Section 7 we give various examples to show that without the condition $\dim \mathcal{G}_K = \dim \mathcal{G}_{\bar{\kappa}}$ imposed in the definition of quasi-reductive group schemes, Theorem 1.2 would be false. In the beginning of Section 9, we give examples of good quasi-reductive models (for definition of a good quasi-reductive model, see Section 8) of $SO_{2n+1}$, $n \geqslant 1$ which are not smooth. Sections 8–10 are devoted to the proof of the fact that our examples are essentially the only such models. In Section 10, using a quadratic form provided to us by Parimala, we construct an example of quasi-reductive model which is not good.

A result of Michel Raynaud (Proposition 3.4) plays a crucial role in this paper. Brian Conrad kindly provided us a proof of it based on an argument of Faltings. With his permission, we have reproduced this proof in the appendix at the end of this paper for the convenience of the reader.

*Acknowledgment.* We thank Kari Vilonen for his question and Ching-Li Chai, Ofer Gabber, Philippe Gille, Michel Raynaud, Jean-Pierre Serre, Marie-France Vigneras, and Shou-Wu Zhang for useful conversations, correspondence and comments. We thank Raman Parimala for the quadratic form used in 10.7 (ii). We are grateful to Brian Conrad for carefully reading several versions of this paper, for his numerous corrections and helpful remarks, and for providing us a proof of Proposition 3.4.

## 2 Unipotent isogenies

Let $k$ be an algebraically closed field of characteristic $p \geqslant 0$. Let $G, G'$ be connected reductive algebraic $k$-groups, and $\phi : G \to G'$ an isogeny. We say that $\phi$ is *unipotent* if the only subgroup scheme of $\ker(\phi)$ of multiplicative type is the trivial subgroup, or, equivalently, $\phi|_T : T \to \phi(T)$ is an isomorphism for a maximal torus $T$.

**2.1 Example** Assume $p = 2$. For $n \geqslant 1$, let $G = SO(V, q) \simeq SO_{2n+1}$, where $(V, q)$ is a quadratic space over $k$ of dimension $2n+1$ and defect 1. Let $\langle \cdot, \cdot \rangle$ be the associated symmetric bilinear form. Then $\langle \cdot, \cdot \rangle$ is also alternating and $V^\perp$ is 1-dimensional. Let $G' = Sp(V/V^\perp) \simeq Sp_{2n}$. Then the natural morphism $\phi_n : G \to G'$ is a unipotent isogeny. The kernel of $\phi_n$ is a finite unipotent group scheme of rank $2^{2n}$.

We will now classify the unipotent isogenies. Let $\phi : G \to G'$ be a unipotent isogeny. Then $\phi$ induces an isomorphism from the connected center $\mathcal{Z}(G)^\circ$ to $\mathcal{Z}(G')^\circ$, and a unipotent isogeny from the derived group $\mathcal{D}(G)$ to $\mathcal{D}(G')$. Let $G_i$, $i \in I$, be the connected normal almost simple algebraic subgroups of $\mathcal{D}(G)$, and $G'_i = \phi(G_i)$, then $G'_i, i \in I$, are the connected normal almost simple algebraic subgroups of $G'$, and $\phi|_{G_i} : G_i \to G'_i$ is a unipotent isogeny for each $i \in I$. Thus it is enough to classify unipotent isogenies between almost simple algebraic groups.

**2.2 Lemma** *Let $\phi : G \to G'$ be a unipotent isogeny between connected almost simple algebraic k-groups, such that $\phi$ is not an isomorphism. Then $p = 2$ and there exists $n \geqslant 1$ such that the morphism $G \xrightarrow{\phi} G'$ is isomorphic to the morphism $\mathrm{SO}_{2n+1} \xrightarrow{\phi_n} \mathrm{Sp}_{2n}$ in the preceding example. That is, there exist isomorphisms $f : G \xrightarrow{\sim} \mathrm{SO}_{2n+1}$, $f' : G' \xrightarrow{\sim} \mathrm{Sp}_{2n}$ such that $f' \circ \phi = \phi_n \circ f$.*

PROOF. Clearly we must have $p > 0$. Let $T$ be a maximal torus of $G$ and $(X, \Phi, X^\vee, \Phi^\vee)$ be the root datum of $(G, T)$, i.e. $X = X^*(T)$ etc. Then we can identify $X$ with $X^*(T')$, where $T' = \phi(T)$. Let $(X, \Phi', X^\vee, \Phi'^\vee)$ be the root datum of $(G', T')$. Then there is a $p$-morphism $\phi_*$ from $(X, \Phi, X^\vee, \Phi^\vee)$ to $(X, \Phi', X^\vee, \Phi'^\vee)$ induced by $\phi$ [SGA3, Exp. XXI, 6.8].

If $\phi$ is not a special isogeny [BoT], then $p^{-1}\Phi' \subset X$. By looking at the classification, we see that this only happens when $p = 2$ and $G \simeq \mathrm{SO}_3$. In this case $G \xrightarrow{\phi} G'$ is isomorphic to $\mathrm{SO}_3 \xrightarrow{\phi_1} \mathrm{Sp}_2$.

If $\phi$ is a special isogeny, then by the classification of special isogenies [BoT], either $p = 3$ and $G$ is of type $G_2$, or $p = 2$ and $G$ is an adjoint group of type $F_4$ or $B_n, n \geqslant 2$. In the case of type $B_n$, $G \xrightarrow{\phi} G'$ is isomorphic to $\mathrm{SO}_{2n+1} \xrightarrow{\phi_n} \mathrm{Sp}_{2n}$.

We claim that $G$ cannot be simply connected, hence we can rule out the cases of type $G_2$ or $F_4$. Indeed, there exists $a \in \Phi$ such that $\ker(\phi|_{U_a})$ is non-trivial, where $U_a$ is the corresponding root subgroup. Let $H$ be the algebraic subgroup of $G$ generated by $U_a$ and $U_{-a}$. Then $\phi$ induces a unipotent isogeny from $H$ to $\phi(H)$. Since $G$ is simply connected, $H$ is isomorphic to $\mathrm{SL}_2$. But every unipotent isogeny from $\mathrm{SL}_2$ is an isomorphism. A contradiction. ■

Note that the root systems of $\mathrm{Sp}_{2n}$ and $\mathrm{SO}_{2n+1}$ are different for $n \geqslant 3$. When $n = 1$ or $2$, the root systems of $\mathrm{Sp}_{2n}$ and $\mathrm{SO}_{2n+1}$ are the same, but $\mathrm{Sp}_{2n}$ is simply connected, whereas $\mathrm{SO}_{2n+1}$ is of adjoint type with nontrivial fundamental group. Therefore, $\mathrm{Sp}_{2n}$ and $\mathrm{SO}_{2n+1}$ are always of different *type* (in the sense discussed after Theorem 1.2).

**2.3 Corollary** *If $\phi : G \to G'$ is a unipotent isogeny which is not an isomorphism, then the type of $G$ is different from that of $G'$.*

**2.4 Corollary** *Let $\phi : G \to G'$ be a unipotent isogeny and $S$ be a maximal torus of $G$. Let $S' = \phi(S)$. If $a$ is an element of $\Phi(G, S)$, the set of roots of $G$ with respect to $S$, then either $a$ or $2a$ is in $\Phi(G', S')$ under the identification $X^*(S) = X^*(S')$.*

## 3 Models

We recall that $R$ is a DVR with residue field $\kappa$. Let $R^{\mathrm{sh}}$ be a strict henselization of $R$. The residue field $\kappa^{\mathrm{sh}}$ of $R^{\mathrm{sh}}$ is then a separable closure of $\kappa$.

Let $G$ be a connected linear algebraic $K$-group. By a *model* of $G$, we mean an affine flat $R$-group scheme $\mathcal{G}$ of finite type over $R$ such that $\mathcal{G}_K := \mathcal{G} \otimes_R K \simeq G$. We have imposed the condition that $\mathcal{G}$ is affine for clarity. According to the following result of Michel Raynaud [SGA3, Exp. XVII, Appendice III, Proposition 2.1(iii)], this is equivalent to $\mathcal{G}$ being separated. Since no published proof of this result is available, for the reader's convenience we are reproducing a proof which was kindly supplied to us by Michel Raynaud.

**3.1 Proposition** *Let $\mathcal{G}$ be a flat group scheme of finite type over $R$ such that its generic fiber $\mathcal{G}_K$ is affine. Then $\mathcal{G}$ is affine if and only if it is separated over $R$.*

PROOF. The "only if" part is obvious. To prove the "if" part, we may replace $\mathcal{G}$ by $\mathcal{G} \otimes_R R'$ for any faithfully flat local extension $R \subset R'$ of DVRs [EGA, IV, Proposition 2.7.1]. Therefore, by Corollary A.4, we may and do assume that $R$ is strictly henselian, and the normalization $\widetilde{\mathcal{G}}_{\text{red}}$ of $\mathcal{G}_{\text{red}}$ is finite over $\mathcal{G}_{\text{red}}$, smooth over $R$, and an $R$-group scheme. Notice that this is where we use the hypothesis that $\mathcal{G}$ is separated over $R$. By a theorem of Chevalley [EGA, II, Proposition 6.7.1], it suffices to show that $\widetilde{\mathcal{G}}_{\text{red}}$ is affine. Therefore, we may and do assume that $\mathcal{G}$ is smooth over $R$.

Let $A = \Gamma(\mathcal{G}, \mathcal{O}_\mathcal{G})$. We remark that it is not clear a priori that $A$ is of finite type over $R$. However, $\mathcal{G}' := \operatorname{Spec} A$ has a natural structure of group scheme over $R$ and the canonical morphism $u : \mathcal{G} \to \mathcal{G}'$ is a morphism of group schemes over $R$. Obviously, $u_K$ is an isomorphism. For $f \in A$, there is a canonical morphism of schemes $u_f : \mathcal{G}_f \to \mathcal{G}'_f$. By [BLR, page 161, Lemma 6], there exists $f \in A$ such that $u_f$ is an isomorphism and $\mathcal{G}_f \cap \mathcal{G}_\kappa \neq \emptyset$. Since translations of $\mathcal{G}_f$ by elements of $\mathcal{G}(R)$ (which maps surjectively to $\mathcal{G}(\kappa)$), together with $\mathcal{G}_K$, form an open covering of $\mathcal{G}$, $u$ is an open immersion. It therefore remains to prove the topological assertion that $u$ is surjective.

By [SGA3, Exp. VI$_B$, Proposition 1.2], $u_\kappa$ is a closed immersion. By [BLR, page 161, before Lemma 6], $A_\kappa \to \Gamma(\mathcal{G}_\kappa, \mathcal{O}_{\mathcal{G}_\kappa})$ is injective. Therefore, $u_\kappa$ is an isomorphism, and hence $u$ is also an isomorphism. ∎

**3.2 Smoothening** For a model $\mathcal{G}$ of $G$, there exists a canonical smoothening morphism $\phi : \widehat{\mathcal{G}} \to \mathcal{G}$ [BLR, 7.1, Theorem 5], which is characterized by the following properties:

(i) $\widehat{\mathcal{G}}$ is a model of $G$, and is smooth over $R$;
(ii) $\widehat{\mathcal{G}} \to \mathcal{G}$ induces a bijection $\widehat{\mathcal{G}}(R^{\text{sh}}) \simeq \mathcal{G}(R^{\text{sh}})$.

In fact, (i) can be replaced by

(i') $\widehat{\mathcal{G}}$ is an affine smooth scheme over $R$ with generic fiber $G$.

Indeed, (i'), (ii), and [BT2, 1.7] show that the group law on the generic fiber of $\widehat{\mathcal{G}}$ extends to $\widehat{\mathcal{G}}$.

We often regard $\mathcal{G}(R)$ as a subgroup of $G(K)$ via the canonical embedding $\mathcal{G}(R) \hookrightarrow G(K)$. For example, (ii) is then simply $\widehat{\mathcal{G}}(R^{\text{sh}}) = \mathcal{G}(R^{\text{sh}})$.

We refer to [BLR] for further properties of the smoothening. For example, [BLR, 7.1, Theorem 5, and 3.6, Proposition 4] imply that the formation of group smoothening is compatible with the base change Spec $R' \to$ Spec $R$ if $R \subset R'$ is a local extension of DVRs of ramification index 1 and the residue field extension of $R'/R$ is separable, such as $R' = R^{\mathrm{sh}}$ (cf. [BLR, 7.2, Theorem 1]).

**3.3 Normalization** Let $\mathcal{G}$ be a model of $G$. The normalization $\widetilde{\mathcal{G}}$ of $\mathcal{G}$ is also an affine flat scheme over $R$ with generic fiber $G$ such that $\widetilde{\mathcal{G}}(R^{\mathrm{sh}}) = \mathcal{G}(R^{\mathrm{sh}})$. Since $\mathcal{G}_K = G$ is an algebraic group, it is geometrically reduced and by Theorem A.6 of the appendix, $\widetilde{\mathcal{G}} \to \mathcal{G}$ *is a finite morphism.*

By the universal property of normalization, the smoothening morphism $\widehat{\mathcal{G}} \to \mathcal{G}$ factors through $\widetilde{\mathcal{G}} \to \mathcal{G}$ uniquely. By 3.2, the morphism $\widehat{\mathcal{G}} \to \widetilde{\mathcal{G}}$ is an isomorphism of schemes if and only if $\widetilde{\mathcal{G}}$ is smooth over $R$. If $\widehat{\mathcal{G}} \simeq \widetilde{\mathcal{G}}$ and $\kappa$ is perfect, *the homomorphism* $\mathcal{G}(R^{\mathrm{sh}}) \to \mathcal{G}(\kappa^{\mathrm{sh}})$ *is surjective.* Indeed, $\widehat{\mathcal{G}}(R^{\mathrm{sh}}) \to \widehat{\mathcal{G}}(\kappa^{\mathrm{sh}})$ is always surjective ($R^{\mathrm{sh}}$ being henselian and $\widehat{\mathcal{G}}$ smooth), and so is $\widetilde{\mathcal{G}}(\kappa^{\mathrm{sh}}) \to \mathcal{G}(\kappa^{\mathrm{sh}})$ (going-up theorem).

It is easy to see that $\widehat{\mathcal{G}} \to \mathcal{G}$ is a finite morphism if and only if $\widehat{\mathcal{G}} \simeq \widetilde{\mathcal{G}}$. Observe that if $\mathcal{G}' \to \mathcal{G}$ is a finite morphism of models of $G$ extending the identity morphism on the generic fibers, then $\mathcal{G}'(R) = \mathcal{G}(R)$. This shows that the condition $\widehat{\mathcal{G}} \simeq \widetilde{\mathcal{G}}$, or equivalently that $\widetilde{\mathcal{G}}$ is smooth over $R$, is stable under the base change Spec $R' \to$ Spec $R$, if $R \subset R'$ is a local extension of DVRs, and in this case the formation of smoothening is stable under such base changes.

**3.4 Proposition** *Assume that $R$ is henselian. Then there exists a local extension $R \subset R'$ of DVRs such that the normalization $\widetilde{\mathcal{G}}'$ of $\mathcal{G}' := \mathcal{G} \otimes_R R'$ is smooth.*

This result is due to Michel Raynaud and appeared in [An, Appendice II, Corollary 3]. Brian Conrad pointed out to us its similarity with Faltings' result ([dJ, Lemma 2.13]), and that Raynaud's result can be deduced by modifying Faltings' argument. For the reader's convenience, we include a complete discussion, which was provided to us by Brian Conrad, in the appendix. In particular, the above result is proved there as Corollary A.4 (also see Remark A.1).

**The neutral component** We recall from [SGA3, Exp. VI$_B$, no. 3] that if $\mathcal{G}$ is a *smooth* group scheme over $R$, the union of $\mathcal{G}_K^\circ$ and $\mathcal{G}_\kappa^\circ$ is an open subgroup scheme of $\mathcal{G}$, which is denoted by $\mathcal{G}^\circ$ and is called the *neutral component* of $\mathcal{G}$.

**3.5 Lemma** *If $\mathcal{G}$ is a smooth affine group scheme over $R$, so is $\mathcal{G}^\circ$.*

PROOF. We need only to show that $\mathcal{G}^\circ$ is affine, which is [BT2, Corollaire 2.2.5 (iii)]. Since $\mathcal{G}^\circ$ is an open subscheme of the affine scheme $\mathcal{G}$, it is separated. Therefore, the result also follows from Proposition 3.1.

Alternatively, one can use the following more elementary claim, together with the fact (see [BLR, 2.1]) that the *dilatation* of an affine scheme is affine: *if $\mathcal{G}$ is a smooth group scheme over $R$ with connected generic fiber, then $\mathcal{G}^\circ$ is the dilatation of $\mathcal{G}_\kappa^\circ$ on $\mathcal{G}$.*

To see this, we observe that the above-mentioned dilatation $\mathcal{G}'$ and $\mathcal{G}^\circ$ are both flat over $R$, and are subfunctors of $\mathcal{G}$ as functors on the category of *flat* $R$-schemes. It is enough to show that they are identical subfunctors. Indeed, the subfunctor $\mathcal{G}^\circ$ is described in [SGA3, Exp. VI$_B$, 3.1], and the subfunctor $\mathcal{G}'$ is described in [BLR, 2.1] (cf. 7.2 below), and they are identical. ∎

*3.6 Chevalley schemes* For simplicity, now assume that $R = R^{\mathrm{sh}}$. Then a reductive group $\mathcal{G}$ over $R$ is necessarily split [SGA3, Exp. XXII, Proposition 2.1], and its isomorphism class is uniquely determined by its type (see the discussion after Theorem 1.2). We will say that $\mathcal{G}$ is a *Chevalley scheme* of that type or the *Chevalley model* of its generic fiber.

The Chevalley model $\mathcal{T}$ of a split torus $T$ over $K$ is unique, and is also called the *Néron-Raynaud model*. It is a *torus over $R$* in the sense of [SGA3, Exp. IX, 1.3], and is characterized by (i) $\mathcal{T}$ is smooth over $R$, and (ii) $\mathcal{T}(R)$ is the maximal bounded subgroup of $T(K)$.

The following lemma is due to Iwahori and Matsumoto, and is now part of the Bruhat–Tits theory [BT2, Proposition 4.6.31]. It is the only fact we need from the Bruhat-Tits theory.

**3.7 Lemma** *Let $\mathcal{G}$ be a Chevalley scheme. Then $\mathcal{G}(R)$ is a hyperspecial maximal bounded subgroup of $\mathcal{G}(K)$. Conversely, if $\mathcal{G}$ is a smooth model of a split connected reductive algebraic group $G$ over $K$ such that $\mathcal{G}(R)$ is a hyperspecial maximal bounded subgroup of $G(K)$, then $\mathcal{G}$ is a Chevalley scheme.*

## 4 First steps of the Proof of Theorem 1.2

In this section, we assume that $R = R^{\mathrm{sh}}$ and $\kappa$ is algebraically closed.

**4.1 Lemma** *Let $\mathcal{T}$ be a model of a split torus $T$ over $K$. Suppose that $\mathcal{T}(R)$ is the maximal bounded subgroup of $T(K)$. Then $\mathcal{T}$ is smooth over $R$, and hence is the Néron-Raynaud model of $T$.*

PROOF. Let $\phi : \widehat{\mathcal{T}} \to \mathcal{T}$ be the smoothening morphism. Then $\widehat{\mathcal{T}}(R) \, (= \mathcal{T}(R))$ is the maximal bounded subgroup of $T(K)$, and hence $\widehat{\mathcal{T}}$ is the Néron-Raynaud model of $T$. We want to show that $\phi$ is an isomorphism.

By [SGA3, Exp. IX, Théorème 6.8], $\ker \phi$ is a group of multiplicative type. Since the generic fiber of $\ker \phi$ is trivial, $\ker \phi$ is trivial. Therefore, $\phi$ is a monomorphism, hence a closed immersion by [SGA3, Exp. IX, Corollaire 2.5]. The ideal (sheaf) $\mathcal{I}$ of this closed immersion has generic fiber 0, since $\phi_K$ is an isomorphism, and so $\mathcal{I} = 0$. Thus, $\phi$ is an isomorphism. ∎

**4.2 Lemma** *Let $\mathcal{G}$ be a model of a connected reductive $K$-group $G$. If $\mathcal{G}(R)$ is a hyperspecial subgroup of $G(K)$, then $\mathcal{G}$ is a quasi-reductive $R$-group scheme and the reduction map $\mathcal{G}(R) \to \mathcal{G}(\kappa)$ is surjective.*

PROOF. Let $\widehat{\mathcal{G}}$ be the smoothening of $\mathcal{G}$. By Lemma 3.7, $\widehat{\mathcal{G}}$ is a Chevalley scheme. For any split maximal torus $\widehat{\mathcal{T}}$ of $\widehat{\mathcal{G}}$, the schematic closure $\mathcal{T}$ of $T := \widehat{\mathcal{T}}_K$ in $\mathcal{G}$ is isomorphic to $\widehat{\mathcal{T}}$ by Lemma 4.1, because $\mathcal{T}(R) = T(K) \cap \mathcal{G}(R) = T(K) \cap \widehat{\mathcal{G}}(R) = \widehat{\mathcal{T}}(R)$. This implies that the kernel of the homomorphism $\widehat{\mathcal{G}}_\kappa \to (\mathcal{G}_\kappa)_{\mathrm{red}}$ does not contain any nontrivial torus and hence the homomorphism is surjective, which implies that $(\mathcal{G}_\kappa)^\circ_{\mathrm{red}}$ is a reductive group. Therefore, $\mathcal{G}$ is a quasi-reductive group scheme. ∎

**Consequences of $\mathcal{G}(R) \twoheadrightarrow \mathcal{G}(\kappa)$**

In the rest of this section, $\mathcal{G}$ is a quasi-reductive model of $G$, $\phi : \widehat{\mathcal{G}} \to \mathcal{G}$ is the smoothening morphism of 3.2, and $\widetilde{\mathcal{G}}$ is the normalization of $\mathcal{G}$.

**4.3 Proposition** *Assume that $\mathcal{G}(R) \to \mathcal{G}(\kappa)$ is surjective. Then*

   (i) *$G$ is a $K$-split reductive group;*
  (ii) *$\mathcal{G}(R)$ is a hyperspecial maximal bounded subgroup of $G(K)$;*
 (iii) *the smoothening $\widehat{\mathcal{G}}$ of $\mathcal{G}$ is a Chevalley scheme;*
 (iv) *there exists an $R$-torus $\mathcal{T}$ in $\mathcal{G}$ such that $\mathcal{T}_K$ is a maximal $K$-split torus of $G$;*
  (v) *the morphism $\widehat{\mathcal{G}}_\kappa \to (\mathcal{G}_\kappa)_{\mathrm{red}}^\circ$ is a unipotent isogeny;*
 (vi) *$G$ is almost simple if and only if $(\mathcal{G}_\kappa)_{\mathrm{red}}^\circ$ is almost simple.*

PROOF. Since the composition $\mathcal{G}(R) = \widehat{\mathcal{G}}(R) \to \widehat{\mathcal{G}}(\kappa) \to \mathcal{G}(\kappa)$ is surjective, $\widehat{\mathcal{G}}(\kappa) \to \mathcal{G}(\kappa)$ is surjective, and so is $(\widehat{\mathcal{G}}_\kappa)^\circ \to (\mathcal{G}_\kappa)_{\mathrm{red}}^\circ$. Since $(\widehat{\mathcal{G}}_\kappa)^\circ$ and $(\mathcal{G}_\kappa)_{\mathrm{red}}^\circ$ have the same dimension, it follows that $(\widehat{\mathcal{G}}_\kappa)^\circ$ is a reductive group, for otherwise, its unipotent radical would be a (connected normal) unipotent subgroup of positive dimension lying in the kernel of the homomorphism, and then the dimension of the image would be strictly smaller. By [SGA3, Exp. XIX, Théorème 2.5], $G$ is a reductive group. Lemma 3.5 implies that $(\widehat{\mathcal{G}})^\circ$ is a reductive group scheme, hence a Chevalley scheme. In particular, $G$ is split over $K$. This proves (i).

Now we have $(\widehat{\mathcal{G}})^\circ(R) \subset \widehat{\mathcal{G}}(R) = \mathcal{G}(R)$ and $(\widehat{\mathcal{G}})^\circ(R)$ is a hyperspecial maximal bounded subgroup of $G(K)$ by Lemma 3.7. So we must have $(\widehat{\mathcal{G}})^\circ(R) = \widehat{\mathcal{G}}(R)$. According to [BT2, 1.7.3c], a smooth model of $G$ is completely determined by its set of integral points. Therefore, we have $(\widehat{\mathcal{G}})^\circ = \widehat{\mathcal{G}}$ and assertions (ii) and (iii) hold.

Since $\widehat{\mathcal{G}}$ is a Chevalley scheme, there exists an $R$-torus $\widehat{\mathcal{T}}$ in $\widehat{\mathcal{G}}$ such that $T := \widehat{\mathcal{T}}_K$ is a maximal $K$-split torus of $G$. Let $\mathcal{T}$ be the schematic closure in $\mathcal{G}$ of $T$. Then $\mathcal{T}(R) = \mathcal{G}(R) \cap \mathcal{T}(K) = \widehat{\mathcal{G}}(R) \cap \widehat{\mathcal{T}}(K) = \widehat{\mathcal{T}}(R)$. By Lemma 4.1, the morphism $\phi|_{\widehat{\mathcal{T}}} : \widehat{\mathcal{T}} \to \mathcal{T}$ is an isomorphism. Thus $\mathcal{T}$ is an $R$-torus of $\mathcal{G}$. This proves (iv).

The morphism between the special fibers induced by $\phi|_{\widehat{\mathcal{T}}}$ is also an isomorphism which implies (v).
Finally, we prove (vi): $G$ is almost simple $\iff \widehat{\mathcal{G}}_\kappa$ is almost simple $\iff (\mathcal{G}_\kappa)_{\mathrm{red}}^\circ$ is almost simple. ∎

**4.4 Proposition** *Let $\mathcal{G}$ be a quasi-reductive model of $G$.*

   (i) *The generic fiber $G$ is a reductive group and $\mathcal{G}_\kappa$ is connected. Moreover, if $H$ is a connected normal $K$-subgroup of $G$, then the schematic closure $\mathcal{H}$ of $H$ in $\mathcal{G}$ is quasi-reductive and $(\mathcal{H}_\kappa)_{\mathrm{red}}$ is a connected normal algebraic subgroup of $(\mathcal{G}_\kappa)_{\mathrm{red}}$.*
  (ii) *Assume that $G$ is $K$-split. Then the correspondence $H \mapsto (\mathcal{H}_\kappa)_{\mathrm{red}}$ is an inclusion preserving bijection from $\mathcal{N}(G)$, the set of connected normal algebraic subgroups of $G$, onto $\mathcal{N}((\mathcal{G}_\kappa)_{\mathrm{red}})$, where for a connected normal algebraic subgroup $H$ of $G$, $\mathcal{H}$ denotes its schematic closure in $\mathcal{G}$.*

PROOF. (i) Let $H$ be a connected normal $K$-subgroup of $G$ and $\mathcal{H}$ be its schematic closure in $\mathcal{G}$. By Proposition 3.4, there exists a local extension $R \subset R'$ of DVRs such that the normalization $\widetilde{\mathcal{G}}'$ of $\mathcal{G}' := \mathcal{G} \otimes_R R'$ and the normalization $\widetilde{\mathcal{H}}'$ of $\mathcal{H}' := \mathcal{H} \otimes_R R'$ are smooth. Then $\mathcal{G}(R') \to \mathcal{G}(\kappa)$ and $\mathcal{H}(R') \to \mathcal{H}(\kappa)$ are surjective (cf. 3.3). We conclude from Proposition 4.3 that if $K'$ is the field of fractions of $R'$, $\mathcal{G}'_{K'} = G \otimes_K K'$, and hence $G$, and so also its connected normal algebraic subgroup $H$, are reductive, and the smoothening $\widehat{\mathcal{G}'}$ ($\simeq \widetilde{\mathcal{G}'}$) is a Chevalley scheme over $R'$. In particular, $\widetilde{\mathcal{G}'}_\kappa$ is connected. As $\mathcal{G}_\kappa = \mathcal{G}'_\kappa$ and $\widetilde{\mathcal{G}'}(\kappa) \to \mathcal{G}'(\kappa)$ is surjective, it follows that $\mathcal{G}_\kappa$ is also connected.

As $\mathcal{H}(R')$ is clearly a normal subgroup of $\mathcal{G}(R')$, $\mathcal{H}(\kappa)$ is a normal subgroup of $\mathcal{G}(\kappa)$, and hence $(\mathcal{H}_\kappa)^\circ_{\mathrm{red}}$ is a normal subgroup of the connected reductive group $(\mathcal{G}_\kappa)_{\mathrm{red}}$ since $\mathcal{G}(\kappa)$ is Zariski-dense in $(\mathcal{G}_\kappa)_{\mathrm{red}}$. Therefore, $(\mathcal{H}_\kappa)^\circ_{\mathrm{red}}$ is reductive and $\mathcal{H}$ is quasi-reductive. This fact now implies that $\mathcal{H}_\kappa$ is actually connected.

(ii) We now assume that $G$ is $K$-split. Then any connected normal algebraic subgroup of $G$ is a reductive group defined over $K$, $H \mapsto H \otimes_K K'$ is a bijection $\mathcal{N}(G) \simeq \mathcal{N}(G \otimes_K K')$ for any field extension $K'/K$. Therefore, to prove (ii) we are free to replace $R$ by a totally ramified local extension $R \subset R'$ of DVRs. Thus we may and do assume that $\widehat{\mathcal{G}}$ is a Chevalley scheme, thanks to Proposition 3.4. Now it is clear from the relation between the root data of $G$, of $\widehat{\mathcal{G}}_\kappa$, and of $(\mathcal{G}_\kappa)_{\mathrm{red}}$, with respect to suitable maximal split tori, that $H \mapsto (\mathcal{H}_\kappa)_{\mathrm{red}}$ is an inclusion preserving bijection from $\mathcal{N}(G)$ to $\mathcal{N}((\mathcal{G}_\kappa)_{\mathrm{red}})$. ∎

**4.5 Proposition** *Let $\mathcal{G}$ be a quasi-reductive model of $G$. Assume that at least one of the three conditions (i), (ii), (iii) of Theorem 1.2 holds. Then $\mathcal{G}$ is smooth over $R$, and is a Chevalley scheme.*

PROOF. It is clear that the hypotheses and the conclusion of the proposition are unchanged if we replace $R$ by $R'$ for any local extension $R \subset R'$ of DVRs. Thanks to Proposition 3.4, by changing $R$, we may assume that the normalization $\widetilde{\mathcal{G}}$ of $\mathcal{G}$ is smooth.

By Proposition 4.3 (v), the isogeny $\widehat{\mathcal{G}}_\kappa \to (\mathcal{G}_\kappa)_{\mathrm{red}}$ is a unipotent isogeny. Therefore, it is an isomorphism by our discussion of unipotent isogenies (§2). Thus, $\phi_\kappa$ is a monomorphism, hence a closed immersion.

Let $A$ and $\widehat{A}$ be the affine rings of $\mathcal{G}$ and $\widehat{\mathcal{G}}$ respectively, and let $\phi^* : A \to \widehat{A}$ be the injective morphism between affine rings. Then $C = \mathrm{coker}(\phi^*)$ is a torsion $R$-module. Since $\phi_\kappa$ is a closed immersion, $C \otimes \kappa = 0$ and hence $C$ is a divisible $R$-module.

But $\widehat{\mathcal{G}} \simeq \widetilde{\mathcal{G}} \to \mathcal{G}$ is a finite morphism (see 3.3). That is, $\widehat{A}$ is a finite $A$-module. If $x_1, \ldots, x_n$ generate $C$ as an $A$-module, and $N \in \mathbb{Z}_{>0}$ is such that $\pi^N x_i = 0$ for all $i$, then $\pi^N C = 0$. This, together with the fact that $C$ is a divisible $R$-module, implies that $C = 0$ and $\widehat{A} = A$. The proposition is proved. ∎

## 5 Proof of Theorem 1.2 (a)

It is clear that the hypotheses and the conclusion of Theorem 1.2 are unchanged if we replace $R$ by $R'$ for any local extension $R \subset R'$ of DVRs. So we may and do assume that $R = R^{\mathrm{sh}}$ and $\kappa$ is algebraically closed. Thus we we have already proved a large part of Theorem 1.2. Indeed, if $\mathcal{G}$ is a quasi-reductive group scheme *of finite type* over $R$, then assertions (b) and (c) of Theorem 1.2 follow from Proposition 4.4 (i), and the last assertion of Theorem 1.2 is Proposition 4.5.

To complete the proof of Theorem 1.2, it remains only to prove assertion (a), i. e., a quasi-reductive group scheme over $R$ is always of finite type over $R$. In this section, we will prove this with the additional assumption that at least one of the three conditions (i), (ii), (iii) of Theorem 1.2 holds. The general case is similar but we need to replace the classification of reductive group schemes with that of good quasi-reductive group schemes (see section 8 for the definition and the result), and is given in 9.7.

**5.1 Lemma** *Let $\mathcal{X}$ be a flat affine scheme over $R$ with affine ring $A$ and generic fiber $X$. We assume that $X$ is an irreducible variety (i.e. an irreducible geometrically reduced $K$-scheme of finite type) and either*

  (i) *$\mathcal{X}(R)$ is Zariski-dense in the generic fiber $X$; or*
  (ii) *there is an ideal $I$ of $A_\kappa := A/\pi A$ such that $A_\kappa/I$ is an integral domain, and*

$$\text{tr. deg}_\kappa(A_\kappa/I) \geqslant \dim X.$$

*Then $A$ contains no non-zero $R$-divisible elements.*

PROOF. Let $J$ be the set of $R$-divisible elements in $A$. Then $J$ is an ideal of $A$, and of $K[X]$. Assume that $J \neq \{0\}$. Then we cannot have (i) since $\mathcal{X}(R)$ is contained in the closed subset defined by $J$.

Now assume (ii). Notice $\dim(K[X]/J) < \dim K[X] = \dim X$ and $(A/J)_\kappa \simeq A_\kappa$. Therefore, we can find $d = \dim X$ elements $\bar{x}_1, \ldots, \bar{x}_d$ in $(A/J)_\kappa/I$ that are algebraically independent over $\kappa$. Lift these elements to $x_1, \ldots, x_d \in A/J$. It is easy to see that $x_1, \ldots, x_d$ are algebraically independent over $K$, this contradicts the fact $\dim(K[X]/J) < d$. ∎

**5.2 Lemma** *Retain the hypothesis of the preceding lemma. Assume furthermore that $R$ is complete. Then $A$ is a free $R$-module.*

PROOF. This follows from the preceding lemma, and the following assertion:

> *Let $V$ be a vector space of at most countable dimension over $K$ and $L$ an $R$-submodule of $V$ such that $L$ contains no non-zero $R$-divisible elements. Then $L$ is a free $R$-module.*

We have been told by Jean-Pierre Serre and Marie-France Vigneras that this assertion is well-known and it appears in an exercise in Bourbaki's *Algebra*, Chap. VII. For the convenience of the reader we give a proof here. We first remark that the completeness of $R$ implies that $\text{Hom}_R(K, K) \simeq \text{Hom}_R(K, K/R)$ and $\text{Ext}^1_R(K, R) = 0$.

We may and do assume that $L \otimes_R K = V$. Let $d = \dim_K V \leqslant \aleph_0$. Let $\{v_i\}_{0 \leqslant i < d}$ be a basis of $V$ over $K$, and put $V_n = K\langle v_0, \ldots, v_{n-1}\rangle$, $L_n = L \cap V_n$. Then $L_1$ is obviously free over $R$. Let $u_0$ be a generator of $L_1$ over $R$. We will construct $\{u_i\}_{0 \leqslant i < d}$ inductively so that $\{u_0, \ldots, u_{n-1}\}$ is a basis of $L_n$ over $R$, as follows.

Assume that we have constructed $u_0, \ldots, u_{n-1}$, and $n + 1 \leqslant d$. Then $L_{n+1}/L_n$ is a non-zero $R$-submodule of $V_{n+1}/V_n \simeq K$, and is isomorphic to either $R$ or $K$. In the latter case, we have $L_{n+1} \simeq L_n \oplus K$ since $\text{Ext}^1_R(K, L_n) = 0$, contradicting the hypothesis on $L$. Therefore, we have a (split) short exact sequence of free $R$-modules $0 \to L_n \to L_{n+1} \to L_{n+1}/L_n \to 0$, from which we construct $u_n$ easily.

It is then clear that $\{u_i : 0 \leqslant i < d\}$ is a basis of $L$ and hence $L$ is free over $R$. ∎

*5.3* We now prove assertion (a) of Theorem 1.2 with with the additional assumption that at least one of the three conditions (i), (ii), (iii) of Theorem 1.2 holds. We may and do assume that $R$ is complete with algebraically closed residue field. Then by the preceding two lemmas, the affine ring $A$ of $\mathcal{G}$ is a free $R$-module. Let $B$ be a basis of $A$ over $R$ and let $S$ be a finite subset of $B$ such that $K[S] = K[G]$ and the image of $S$ in $(A_\kappa)_{\text{red}}$ generates $(A_\kappa)_{\text{red}}$ as a $\kappa$-algebra.

For any finite subset $I$ of $B$ such that $I \supset S$, by the argument in [Wa, 3.3], there is a Hopf subalgebra $A_I$ of $A$ which is a finitely generated $R$-algebra containing $I$. Let $\mathcal{G}^I = \operatorname{Spec} A_I$ be the affine group scheme with affine ring $A_I$. Then as $S \subset I$, the reduced special fiber $(\mathcal{G}_\kappa^I)_{\text{red}}$ of $\mathcal{G}^I$ contains $(\mathcal{G}_\kappa)_{\text{red}}$ as a closed subgroup; in particular, $\dim(\mathcal{G}_\kappa^I)_{\text{red}} \geqslant \dim(\mathcal{G}_\kappa)_{\text{red}}$. By [M2, 15.3], $\dim(\mathcal{G}_\kappa^I)_{\text{red}} \leqslant \dim G = \dim(\mathcal{G}_\kappa)_{\text{red}}$ and hence, $\dim(\mathcal{G}_\kappa^I)_{\text{red}} = \dim(\mathcal{G}_\kappa)_{\text{red}}$. This implies that $(\mathcal{G}_\kappa^I)_{\text{red}}^\circ = (\mathcal{G}_\kappa)_{\text{red}}^\circ$. Therefore, each $\mathcal{G}^I$ is a quasi-reductive group scheme of finite type over $R$. Now we assume that at least one of the three additional conditions of Theorem 1.2 holds. Then it follows from Proposition 4.5 that $\mathcal{G}^I = \operatorname{Spec} A_I$ is a Chevalley scheme.

For any finite sets $I, J$ such that $S \subset I, J \subset B$, we can find a Hopf subalgebra $A'$ of $A$ which contains both $A_I$ and $A_J$, and which is a finitely generated $R$-algebra. Again, $\mathcal{G}' := \operatorname{Spec} A'$ is a Chevalley scheme and the morphism of Chevalley schemes $\mathcal{G}' \to \mathcal{G}^I$ is an isomorphism by Lemma 3.7. Thus $\mathcal{G}^I = \mathcal{G}^J$ and $A_I = A_J$. Since $A$ is the union of $A_I$ for varying $I$, the theorem is proved. ∎

*5.4 Proof of Corollary 1.3* Again, we may and do assume that $R$ is strictly henselian and $\kappa$ is algebraically closed. Let $\mathcal{G}^*$ be the schematic closure of $\phi(\mathcal{G}_K)$ in $\mathcal{H}$. Then $\mathcal{G}^*$ is a model of $G := \mathcal{G}_K$ and $\phi$ factors as $\mathcal{G} \to \mathcal{G}^* \hookrightarrow \mathcal{H}$.

By Proposition 3.4, we can find a local extension $R \subset R'$ of DVRs such that $\mathcal{G}^*(R') \to \mathcal{G}^*(\kappa)$ is surjective. Let $K'$ be the field of fractions of $R'$. Since $\mathcal{G}(R')$ is a maximal bounded subgroup of $G(K')$ (Lemma 3.7) and $\mathcal{G}(R') \subset \mathcal{G}^*(R') \subset G(K')$, we have $\mathcal{G}(R') = \mathcal{G}^*(R')$ and $\mathcal{G}(\kappa) \to \mathcal{G}^*(\kappa)$ is surjective. It follows that $(\mathcal{G}_\kappa^*)_{\text{red}}^\circ$ is a reductive group. Now Theorem 1.2 and [BT2, 1.7] imply that $\mathcal{G} \to \mathcal{G}^*$ is an isomorphism, and hence, $\mathcal{G} \to \mathcal{G}^*$ is also an isomorphism. ∎

## 6 General noetherian base schemes

We will now give the analogues of Theorem 1.2 over more general base schemes. Proposition 6.1 gives a fiberwise result. Theorem 6.2 is a local result and is the reductive analogue of Theorem 1.1 for group schemes of finite type. This is then globalized into Theorem 6.3. Finally, we give the proof of Proposition 1.5.

For the sake of convenience, we introduce the following condition: let $\mathcal{G}$ be a group scheme over a base scheme $S$. Given $s, t \in S$ such that $(\mathcal{G}_{\bar{s}})_{\text{red}}^\circ$ and $(\mathcal{G}_{\bar{t}})_{\text{red}}^\circ$ are reductive algebraic groups for some (hence any) geometric points $\bar{s}$ (resp. $\bar{t}$) over $s$ (resp. $t$), we say that $(\mathcal{G}, s, t)$ *satisfies condition* $(*)$ if at least one of the following holds:

- the characteristic of $\kappa(s)$ is not 2;
- the type of $(\mathcal{G}_{\bar{t}})_{\text{red}}^\circ$ is the same as that of $(\mathcal{G}_{\bar{s}})_{\text{red}}^\circ$;

- no normal algebraic subgroup of $(\mathcal{G}_{\bar{t}})^\circ_{\mathrm{red}}$ is isomorphic to $\mathrm{SO}_{2n+1}$ for $n \geqslant 1$.

**6.1 Proposition** *Let $\mathcal{G}$ be an affine flat group scheme over a noetherian scheme $S$. Let $s \in S$, $\bar{s}$ a geometric point over $s$, and suppose*

- *$(\mathcal{G}_{\bar{s}})^\circ_{\mathrm{red}}$ is a reductive group; and*
- *there exists a generization $t$ of $s$ such that $\mathcal{G}_t$ is connected, smooth over $\kappa(t)$ of dimension equal to that of $(\mathcal{G}_{\bar{s}})^\circ_{\mathrm{red}}$.*

*Then $\mathcal{G}_t$ is a reductive group over $\kappa(t)$, and $\mathcal{G}_s$ is connected. If, in addition, $(\mathcal{G}, s, t)$ satisfies condition $(*)$, then $\mathcal{G}_s$ is connected and smooth over $\kappa(s)$.*

PROOF. There is nothing to prove if $s = t$. So assume $s \neq t$. By [EGA, II, Proposition 7.1.9], we can find a DVR $R$, and a morphism $\mathrm{Spec}\, R \to S$ sending the generic point of $\mathrm{Spec}\, R$ to $t$, and the special point of $\mathrm{Spec}\, R$ to $s$. Now Theorem 1.2 can be applied to $\mathcal{G} \times_S \mathrm{Spec}\, R$ and the proposition follows immediately. ∎

**6.2 Theorem** *Let $\mathcal{G}$ be an affine flat group scheme of finite type over a noetherian scheme $S$, with connected fibers of the same dimension. Let $s \in S$, $\bar{s}$ a geometric point over $s$, and suppose*

- *$(\mathcal{G}_{\bar{s}})_{\mathrm{red}}$ is a reductive group; and*
- *there exists a generization $t$ of $s$ such that $\mathcal{G}_t$ is smooth over $\kappa(t)$.*

*Then $\mathcal{G}_t$ is a reductive group over $\kappa(t)$. If, in addition, we assume that $(\mathcal{G}, s, t)$ satisfies condition $(*)$, then there exists an open neighborhood $U$ of $s$ such that $\mathcal{G}|U$ is a reductive group over $U$.*

PROOF. That $\mathcal{G}_t$ is a reductive group over $\kappa(t)$ is immediate from the preceding proposition, which also implies that $\mathcal{G}_s$ is smooth over $\kappa(s)$. The $S$-smooth locus in $\mathcal{G}$ is a union of fibers since $\mathcal{G}$ is a flat and finite-type $S$-group, and this smooth locus has open image in $S$ again since $\mathcal{G}$ is flat and of finite type over $S$, so $\mathcal{G}$ is smooth over an open neighborhood $U'$ of $s$ (cf. [SGA3, Exp. X, Lemme 3.5]). Hence $\mathcal{G}$ is a reductive group over a smaller neighborhood $U$ of $s$ by [SGA3, Exp. XIX, Théorème 2.5]. ∎

**6.3 Theorem** *Let $\mathcal{G}$ be an affine flat group scheme over an irreducible noetherian scheme $S$, with fibers of finite type and the same dimension. Let $\xi$ be the generic point of $S$ and assume that the generic fiber $\mathcal{G}_\xi$ is smooth and connected. For any geometric point $\bar{s}$ over a non-generic point $s \in S$, assume that $(\mathcal{G}_{\bar{s}})^\circ_{\mathrm{red}}$ is a reductive group. Then $\mathcal{G}_{\bar{s}}$ is connected for all $s \in S$.*

*If, in addition, we assume that $(\mathcal{G}, s, \xi)$ satisfies condition $(*)$ for all non-generic $s \in S$, then each geometric fiber of $\mathcal{G}$ over $S$ is a connected reductive algebraic group. Furthermore, if $\mathcal{G}$ is of finite type over $S$, then $\mathcal{G}$ is a reductive group over $S$.*

PROOF. The first two assertions follow immediately from Proposition 6.1. The last assertion is clear from Theorem 6.2. ∎

**6.4 Proof of Proposition 1.5** It suffices to show that $\mathcal{G}$ is of finite type over $S$, or equivalently, to show that for each closed point $s \in S$, there is an open neighborhood $U$ of $s$ such that $\mathcal{G}|U$ is of finite type over $U$.

Let $C = \mathcal{O}_s$. By Theorem 1.2, $\mathcal{G} \times_S \operatorname{Spec} C$ is of finite type over $\operatorname{Spec} C$, and hence it is a reductive group over $\operatorname{Spec} C$. We can "spread out" $\mathcal{G} \times_S \operatorname{Spec} C$ to an affine group scheme of finite type $\mathcal{G}'$ over an open neighborhood $U'$ of $s$ such that there is a morphism of group schemes $\mathcal{G}|U' \to \mathcal{G}'$ inducing an isomorphism $\mathcal{G}' \times_S \operatorname{Spec} C \simeq \mathcal{G} \times_S \operatorname{Spec} C$. (Concretely, we may assume that $S = \operatorname{Spec} A$ is affine, and $\mathcal{G} = \operatorname{Spec} B$, where $B$ is a Hopf $A$-algebra. Let $x_1, \ldots, x_n$ be elements generating $B \otimes_A C$ over $C$; then there exists an $f \in A$, $f$ a unit in $C$, such that $B' := A_f[x_1, \ldots, x_n] \subset B \otimes_A C$ is a Hopf subalgebra. We then take $U' = \operatorname{Spec} A_f$, $\mathcal{G}' = \operatorname{Spec} B'$.) We may assume that $\mathcal{G}'$ is smooth over $U'$ by shrinking $U'$.

By [SGA3, Exp. XIX, Théorème 2.5], there is an open neighborhood $U$ of $s$ such that $\mathcal{G}'|U$ is a reductive group. We now claim that the morphism $\mathcal{G}|U \to \mathcal{G}'|U$ is an isomorphism, hence $\mathcal{G}|U$ is of finite type over $U$.

If suffices to show that for each closed point $t$ of $U$, $\mathcal{G} \times_S \operatorname{Spec} R$ is isomorphic to $\mathcal{G}' \times_U \operatorname{Spec} R$, where $R = \mathcal{O}_t$, or equivalently, $\mathcal{G}_{\widetilde{R}} := \mathcal{G} \times_S \operatorname{Spec} \widetilde{R}$ is isomorphic to $\mathcal{G}'_{\widetilde{R}} := \mathcal{G}' \times_U \operatorname{Spec} \widetilde{R}$, where $\widetilde{R}$ is the strict henselization of $R$. But since $\mathcal{G}(\widetilde{R})$ is a (hyperspecial) maximal bounded subgroup (Theorem 1.2 and Lemma 3.7), we must have $\mathcal{G}(\widetilde{R}) = \mathcal{G}'(\widetilde{R})$. As both $\mathcal{G}_{\widetilde{R}}$ and $\mathcal{G}'_{\widetilde{R}}$ are affine smooth over $\widetilde{R}$, we have $\mathcal{G}_{\widetilde{R}} = \mathcal{G}'_{\widetilde{R}}$ by [BT2, 1.7]. The proposition is proved. ∎

## 7 Examples

**7.1** The example in 9.1 and 9.2 shows that the hypotheses of Corollary 1.3 are necessary. For $c \neq 0$, the morphism $\mathcal{G}_0 \to \mathcal{H}$ is not a closed immersion, although it is so at the generic fibers.

We will now give examples to show that if from the definition of quasi-reductive group schemes we drop the condition that $\dim(\mathcal{G}_{\bar{\kappa}})_{\mathrm{red}} = \dim \mathcal{G}_K$, then Theorem 1.2 is false. These examples are constructed using variations of *dilatations* ([BLR, 3.2]).

Throughout this section, $R$ is a DVR and $\pi$ is a uniformizer of $R$.

**7.2 Higher dilatations** Let $\mathcal{X}$ be a flat scheme of finite type over $R$, and $\mathcal{Z} \hookrightarrow \mathcal{X}$ be a flat closed subscheme over $R$.

We define a sequence of flat schemes $\Gamma_n = \Gamma_n(\mathcal{X}, \mathcal{Z})$ over $R$, together with closed immersions $i_n : \mathcal{Z} \hookrightarrow \Gamma_n$ as follows. Let $\Gamma_0(\mathcal{X}, \mathcal{Z}) = \mathcal{X}$, and $i_0 : \mathcal{Z} \hookrightarrow \Gamma_0$ be the inclusion. After $\Gamma_n$ and $i_n$ have been defined, we let $\Gamma_{n+1}$ be the dilatation of $i_n(\mathcal{Z}_\kappa)$ on $\Gamma_n$. The dilatation of $\mathcal{Z}_\kappa$ on $\mathcal{Z}$, which is nothing but $\mathcal{Z}$ itself, then admits a natural closed immersion into $\Gamma_{n+1}$ by [BLR, 3.2, Proposition 2(c)], which we denote by $i_{n+1}$.

The scheme $\Gamma_n$ is flat over $R$ by construction, hence is determined by its associated functor of points on the category of flat $R$-schemes. In fact,

**7.3 Proposition** (i) *For any flat $R$-scheme $\mathcal{Y}$, $\Gamma_n(\mathcal{Y})$, the set of $\mathcal{Y}$-valued points in $\Gamma_n$, is*

$$\left\{\phi: \mathcal{Y} \to \mathcal{X} \;\middle|\; \begin{array}{l} \phi \otimes_R (R/\pi^n R) \text{ factors through} \\ \mathcal{Z} \otimes_R (R/\pi^n R) \hookrightarrow \mathcal{X} \otimes_R (R/\pi^n R) \end{array}\right\}$$

(ii) *By (i), each $\Gamma_n$ is a subfunctor of $\mathcal{X}$, and we have $\Gamma_{n+1} \subset \Gamma_n$ for all $n$. Let $\Gamma_\infty = \bigcap_{n \geq 0} \Gamma_n$. Then the functor $\Gamma_\infty$ is also represented by a flat scheme over $R$.*

(iii) *The special fiber of $\Gamma_\infty$ is naturally isomorphic to $\mathcal{Z}_\kappa$.*

PROOF. The statements are local, and can be checked by assuming that $\mathcal{X} = \operatorname{Spec} A$ is affine. Then $\mathcal{Z} = \operatorname{Spec}(A/I)$ for some ideal $I$ of $A$.

By the definition of dilatations ([BLR, 3.2]), $\Gamma_1$ is affine with affine ring

$$A_1 = A[\pi^{-1}I] = A + \sum_{k \geq 1} \pi^{-k} I^k \subset A \otimes_R K.$$

From this, one verifies inductively that $\Gamma_n$ is affine with affine ring

$$A_n = A[\pi^{-n}I] = A + \sum_{k \geq 1} \pi^{-kn} I^k \subset A \otimes_R K.$$

Now it is easy to verify that (i) describes the functor of points with values in flat $R$-algebras.

It is easy to check that $\Gamma_\infty$ is represented by $\operatorname{Spec} A_\infty$, where

$$A_\infty = \bigcup_{n \geq 1} A_n = A + I \otimes_R K,$$

which proves (ii). Now the exact sequence $0 \to I \otimes_R K \to A_\infty \to B \to 0$ shows that $(A_\infty)_\kappa \simeq B_\kappa$, hence (iii). ∎

**7.4 Group schemes** Now assume that $\mathcal{X}$ is a group scheme over $R$, flat of finite type, and $\mathcal{Z}$ is a flat closed subgroup scheme. Then each $\Gamma_n$ is naturally a group scheme over $R$ as well. Indeed, this follows from [BLR, 3.2, Proposition 2(d)]. One can also use the observation that the subfunctor $\Gamma_n$ of $\mathcal{X}$ is a subgroup functor. When $\mathcal{Z}$ is the trivial subgroup, $\Gamma_n$ is the scheme-theoretic principal congruence subgroup of level $n$. In general, we can think of $\Gamma_n$ as a scheme-theoretic congruence subgroup.

Similarly, $\Gamma_\infty$ is naturally a group scheme over $R$.

**7.5 Examples** Let $\mathcal{H}$ be any flat group scheme of finite type over $R$, and $\mathcal{E}$ be the trivial closed subgroup scheme. Then $\mathcal{G} := \Gamma_\infty(\mathcal{H}, \mathcal{E})$ has the same generic fiber as $\mathcal{H}$, but its special fiber is the trivial group. In particular, the special fiber is reductive.

We can also take $\mathcal{H}$ to be a Chevalley scheme over $R$, and $\mathcal{H}'$ a proper Levi subgroup. Then again we see that the special fiber of $\mathcal{G} := \Gamma_\infty(\mathcal{H}, \mathcal{H}')$ is reductive. These examples show that without the condition that $\dim(\mathcal{G}_{\bar{\kappa}})_{\mathrm{red}} = \dim \mathcal{G}_K$ in the definition of quasi-reductive group schemes, Theorem 1.2 is false.

**7.6 A variant** Let $\mathcal{H}$ be a flat group scheme of finite type over $R$. For simplicity, we assume that $\mathcal{H}$ is affine and $\kappa$ is a finite field with $q$ elements.

Let $F : \mathcal{H}_\kappa \to \mathcal{H}_\kappa$ be the Frobenius $\kappa$-morphism, and $Z = \ker(F)$. We define $\Lambda_1 = \Lambda_1(\mathcal{H})$ to be the dilatation of $Z$ on $\mathcal{H}$. In the following, we use $\Gamma_1(\mathcal{H})$ to denote $\Gamma_1(\mathcal{H}, \mathcal{E})$, where $\mathcal{E}$ is the trivial subgroup scheme, so that $\Gamma_1(\mathcal{H})(R)$ is the first principal congruence subgroup of $\mathcal{H}$.

**Claim 1.** We have $\Gamma_1(\mathcal{H})(R) \subset \Gamma_1(\Lambda_1)(R)$.

PROOF. We recall that an element $h \in \mathcal{H}(R)$ lies in $\Gamma_1(\mathcal{H})(R)$ if and only if for any $f \in A = \mathcal{O}_\mathcal{H}(\mathcal{H})$ such that $f(e) \in \pi R$, we have $f(h) \in \pi R$, where $e$ is the identity of $\mathcal{H}(R)$.

Notice also that $Z$ is defined by the ideal generated by $\pi$ and $\{f^q : f \in A, f(e) \in \pi R\}$. Therefore, the affine ring of $\Lambda_1$ is $A_1 := A[f^q/\pi : f \in A, f(e) \in \pi R]$.

Now assume that $h \in \Gamma_1(\mathcal{H})(R)$ and $f_1 \in A_1$ is such that $f_1(e) \in \pi R$. We need to show that $f_1(h) \in \pi R$. We can write $f_1 = F(g_1^q/\pi, \ldots, g_n^q/\pi)$, where $g_i \in A$ are such that $g_i(e) \in \pi R$ and $F$ is a polynomial over $A$ in $n$ variables. It follows immediately that the constant term $c_0$ of $F$ satisfies $c_0(e) \in \pi R$. Therefore, we have $c_0(h) \in \pi R$ and $g_i(h) \in \pi R$ for all $i$, and hence $f_1(h) \in \pi R$. ∎

Now let $\Lambda_0 := \mathcal{H}$, and define $\Lambda_{n+1} := \Lambda_1(\Lambda_n)$ inductively. Then the affine rings $A_n$ of the schemes $\Lambda_n$ form an increasing sequence of Hopf algebras (inside the affine ring of $\mathcal{H} \otimes_R K$). Their union $A_\infty$ is again a Hopf algebra, whose spectrum is a group scheme $\Lambda_\infty$ over $R$. By Claim 1, $\Gamma_1(\mathcal{H})(R) \subset \Lambda_\infty(R)$.

**Claim 2.** The reduced special fiber of $\Lambda_\infty$ is the trivial group.

PROOF. Let $f \in A_\infty/\pi A_\infty$. It suffices to show that if $f(e) = 0$, then $f$ is nilpotent.

Lift $f$ to $\tilde{f} \in A_\infty$. Then $\tilde{f} \in A_n$ for some $n \geqslant 0$. Since $\tilde{f}(e) \in \pi R$, $\tilde{f}^q/\pi \in A_{n+1}$. Thus $\tilde{f}^q \in \pi A_\infty$ and hence $f^q = 0$. ∎

If we apply this to $\mathcal{H} = \mathbb{G}_a/R$, then $\Lambda_\infty = \operatorname{Spec} A_\infty$, where $A_\infty = R[X_0, X_1, \ldots] \subset K[X]$, with $X_0 = X$, $X_{n+1} = X_n^q/\pi$. Therefore, $\mathcal{G} := \Lambda_\infty$ has reductive reduced special fiber (trivial group), $\mathcal{G}(R)$ is Zariski-dense in $\mathcal{G} \otimes_R K$, but the generic fiber is not reductive.

We can apply this to a Chevalley scheme $\mathcal{H}$ of positive dimension, and put $\mathcal{G} = \Lambda_\infty$. Then we see again that without the condition $\dim(\mathcal{G}_{\bar{\kappa}})_{\mathrm{red}} = \dim \mathcal{G}_K$, in the definition of quasi-reductive group schemes, Theorem 1.2 is false even under the additional hypothesis that $\mathcal{G}_K$ is reductive and $\mathcal{G}(R)$ is Zariski-dense.

## 8 Good quasi-reductive group schemes

In the remaining three sections (8, 9, 10), we assume that $R$ is a strictly henselian DVR with algebraically closed residue field, and $K$ is its field of fractions. Let $G$ be a connected reductive $K$-group.

**8.1 Lemma** *Let $\mathcal{G}$ be a quasi-reductive model of $G$ over $R$. Consider the following conditions:*

(1) $\widehat{\mathcal{G}} \to \widetilde{\mathcal{G}}$ *is an isomorphism.*
(2) $\mathcal{G}(R) \to \mathcal{G}(\kappa)$ *is surjective.*

(3) $\mathcal{G}(R)$ *is a hyperspecial maximal bounded subgroup of* $G(K)$.

(4) *$G$ is $K$-split and there is an $R$-torus $\mathcal{T}$ in $\mathcal{G}$ such that $\mathcal{T}_K$ is a maximal $K$-split torus of $G$.*

*Then* (1) $\Rightarrow$ (2) $\Leftrightarrow$ (3) $\Rightarrow$ (4).

PROOF. This follows from Lemma 4.2 and Proposition 4.3. ∎

**8.2 Definition** A quasi-reductive model $\mathcal{G}$ of $G$ is called *good* if the equivalent conditions (2) and (3) of the above lemma hold. If $G$ admits a good quasi-reductive model, then it splits over $K$.

We will say that a quasi-reductive $R$-group scheme of finite type is good if it is a good quasi-reductive model of its generic fiber.

We will see later (Theorems 9.3 and 9.4) that the three conditions (1)–(3) are, in fact, equivalent. Notice that by Proposition 3.4, for any quasi-reductive model $\mathcal{G}$ of $G$ over $R$, there is a local extension $R \subset R'$ of DVRs such that $\mathcal{G}' := \mathcal{G} \otimes_R R'$ is good.

**The $\Delta$-invariant** For a $R$-group scheme $\mathcal{G}$ of finite type with connected reductive generic fiber $G$, let $\Xi := \Xi(G)$ be the set of normal algebraic subgroups of $G$ which are isomorphic to $SO_{2n+1}$ for some $n \geqslant 1$. For $H \in \Xi$, define

$$\Delta_{\mathcal{G}}(H) := \frac{1}{2n} \operatorname{length}_R(\operatorname{Lie} \mathcal{H}/\operatorname{Lie} \widetilde{\mathcal{H}}),$$

where $\mathcal{H}$ is the schematic closure of $H$ in $\mathcal{G}$, $\widetilde{\mathcal{H}}$ is the normalization of $\mathcal{H}$, and $n$ is such that $H \simeq SO_{2n+1}$.

If $G$ is an almost simple group, $\Delta_{\mathcal{G}}(G)$ is, by definition, the $\Delta$-invariant of $\mathcal{G}$.

**8.3 Theorem** *Let $\mathcal{G}$ and $\mathcal{G}'$ be good quasi-reductive $R$-group schemes of finite type. Then $\mathcal{G}$ is isomorphic to $\mathcal{G}'$ if and only if there is an isomorphism from $\mathcal{G}_K$ to $\mathcal{G}'_K$ inducing a bijection $\xi : \Xi(\mathcal{G}_K) \to \Xi(\mathcal{G}'_K)$ such that $\Delta_{\mathcal{G}} = \Delta_{\mathcal{G}'} \circ \xi$. Moreover, $\Delta$ takes values in*

$$I = \{m \in \mathbb{Z} : 0 \leqslant 2m \leqslant \operatorname{ord}_K(2)\},$$

*and $\mathcal{G}$ is a reductive group over $R$ if and only if $\Delta_{\mathcal{G}}(H) = 0$ for all $H \in \Xi(\mathcal{G}_K)$. For a fixed $K$-split reductive group $G$, the set of isomorphism classes of good quasi-reductive $R$-group schemes with generic fiber $\simeq G$ is in bijection with*

$$I^{\Xi} := \text{the set of all functions } \Xi \to I,$$

*where $\Xi = \Xi(G)$.*

The proof of this theorem will be given in section 9.

**Good tori and the big-cell decomposition** Let $\mathcal{G}$ be a model of a connected reductive $K$-group $G$. We assume that $G$ is $K$-split for simplicity. We define a *good torus* of $\mathcal{G}$ to be a closed subgroup $\mathcal{S} \hookrightarrow \mathcal{G}$ such that $\mathcal{S}$ is a (split) torus over $R$ and $\mathcal{S}_K$ is a maximal $K$-split torus of $G$. To give such a torus is to give a maximal $K$-split torus $S$ of $G$, and to check that the schematic closure $\mathcal{S}$ of $S$ in $\mathcal{G}$ is the Néron-Raynaud model of $S$. Such a torus may or may not exist in a model $\mathcal{G}$. We now assume that $\mathcal{G}$ contains a good torus $\mathcal{S}$. Let $S$ be the generic fiber of this torus.

Let $\Phi = \Phi(G, S)$ be the root system of $G$ with respect to $S$. For $a \in \Phi$, let $U_a$ be the corresponding root subgroup and $\mathcal{U}_a$ the schematic closure of $U_a$ in $\mathcal{G}$.

**8.4 Theorem** *Assume that $\mathcal{S}$ is a good torus of $\mathcal{G}$. With the above notation, for any system of positive roots $\Phi^+$ of $\Phi$, the multiplication morphism*

$$\left(\prod_{a \in \Phi^+} \mathcal{U}_a\right) \times \mathcal{S} \times \left(\prod_{a \in \Phi^-} \mathcal{U}_a\right) \to \mathcal{G}$$

*is an open immersion, here $\Phi^- = -\Phi^+$, and the two products can be taken in any order.*

This follows from [BT2, 1.4.5] and [BT2, Théorème 2.2.3]. This result is the big-cell decomposition for models with a good torus. We remark that no smoothness assumption about $\mathcal{G}$ is needed here, and the result can be formulated for models of linear algebraic groups that are not reductive. For similar results for certain models without a good torus, see [Yu].

**8.5** Let $\mathcal{G}$ be a good quasi-reductive model of $G$. By Lemma 8.1, there exists a good torus $\mathcal{S}$ of $\mathcal{G}$. By Lemma 4.1, $\mathcal{S}$ can be regarded as a good torus of $\widehat{\mathcal{G}}$ as well. Let $S = \mathcal{S}_K$, $\Phi = \Phi(G, S)$. For $a \in \Phi(G, S)$, let $U_a$ be the corresponding root subgroup. As $G$ splits over $K$ (Lemma 8.1), $U_a$ is $K$-isomorphic to $\mathbb{G}_a$. Let $\mathcal{U}_a$ (resp. $\widehat{\mathcal{U}}_a$) be the schematic closure of $U_a$ in $\mathcal{G}$ (resp. $\widehat{\mathcal{G}}$). Since $\widehat{\mathcal{G}}$ is a Chevalley scheme, $\widehat{\mathcal{U}}_a$ is smooth over $R$.

Since $\widehat{\mathcal{U}}_a \simeq \mathbb{G}_a/R$, the affine ring $R[\widehat{\mathcal{U}}_a]$ is a polynomial ring $R[x]$, we may and do assume that $R[x] = \bigoplus_{i \geq 0} Rx^i$ is the weight decomposition of the action of $\mathcal{S}$ on $R[\widehat{\mathcal{U}}_a]$, where the weight of $Rx^i$ is $ia$ [SGA3, Exp. I, 4.7.3].

The affine ring $R[\mathcal{U}_a]$ is a subring of $R[\widehat{\mathcal{U}}_a]$ and it has a similar weight decomposition $R[\mathcal{U}_a] = \bigoplus_{i \geq 0} R[\mathcal{U}_a] \cap Rx^i$.

**8.6 Proposition** *Let $R[\mathcal{U}_a] \cap Rx^i = R\pi^{n_i} x^i$. Then the following assertions hold.*

(i) *There exists $i \geq 1$ such that $n_i = 0$.*
(ii) *The weights of $\mathrm{Lie}(\mathcal{U}_a \otimes_R \kappa)$ for the action of $\mathcal{S}_\kappa$ are of the form $ia$ with $i \geq 1$.*
(iii) *Let $j = \min\{i \geq 1 : n_i = 0\}$. The weight of the 1-dimensional Lie algebra $\mathrm{Lie}(\mathcal{U}_a \otimes_R \kappa)_{\mathrm{red}}$, for the action of $\mathcal{S}_\kappa$, is $ja$.*
(iv) *If $a$ is the weight of $\mathrm{Lie}(\mathcal{U}_a \otimes_R \kappa)_{\mathrm{red}}$, then $\mathcal{U}_a \simeq \widehat{\mathcal{U}}_a$ is smooth.*

(v) *If $2a$ is the weight of $\mathrm{Lie}(\mathcal{U}_a \otimes_R \kappa)_{\mathrm{red}}$, then there exists an integer $c$ such that $2 \leqslant 2c \leqslant \mathrm{ord}_K(2)$ and $R[\mathcal{U}_a] = R[\pi^c x, x^2]$. In this case,*

$$\mathrm{length}_R(\mathrm{Lie}\,\mathcal{U}_a/\mathrm{Lie}\,\widehat{\mathcal{U}}_a) = c.$$

PROOF. Put $y_i = \pi^{n_i} x^i$. We may assume that $R[\mathcal{U}_a]$ is generated by $y_i$ for $i = 1, \ldots, m$.

(i) Assume the contrary, then $n_i \geqslant 1$ for $i = 1, \ldots, m$. Let $R' = R[\pi^{1/m}]$, $K' = \mathrm{Frac}\,R'$. Then the $R$-algebra homomorphism $R[\mathcal{U}_a] \to R'$, $x \mapsto \pi^{-1/m}$ (i.e. $y_i \mapsto \pi^{n_i-(i/m)}$) is well-defined. This gives a point of $\mathcal{U}_a(R')$ which is not in $\widehat{\mathcal{U}}_a(R')$. However, this contradicts the following observation.

For any totally ramified local extension $R \subset R'$ of DVRs, it is clear from condition (2) of Lemma 8.1 that $\mathcal{G} \otimes_R R'$ is good. Therefore, $\mathcal{G}(R')$ is a hyperspecial maximal bounded subgroup. Since $\widehat{\mathcal{G}}(R') \subset \mathcal{G}(R')$ and $\widehat{\mathcal{G}}$ is a Chevalley group, we must have $\widehat{\mathcal{G}}(R') = \mathcal{G}(R')$. Consequently, $\widehat{\mathcal{U}}_a(R') = \widehat{\mathcal{G}}(R') \cap U_a(K') = \mathcal{G}(R') \cap U_a(K') = \mathcal{U}_a(R')$, where $K' = \mathrm{Frac}\,R'$.

(ii) By definition, an element of $\mathrm{Lie}(\mathcal{U}_a \otimes_R \kappa)$ is a homomorphism $R[\mathcal{U}_a] \otimes_R \kappa \to \kappa[\epsilon]/(\epsilon^2)$, sending each $y_i \otimes 1$ to an element of $\kappa \cdot \epsilon$. Now (ii) is obvious from this description.

(iii) By (i), if $n_i > 0$, then $y_i \otimes \kappa$ is nilpotent in $R[\mathcal{U}_a] \otimes \kappa$. Therefore, $(R[\mathcal{U}_a] \otimes \kappa)_{\mathrm{red}}$ is spanned by the $y_i \otimes 1$ for all $i$ such that $n_i = 0$. An element of $\mathrm{Lie}((\mathcal{U}_a \otimes_R \kappa)_{\mathrm{red}})$ is a homomorphism $(R[\mathcal{U}_a] \otimes_R \kappa)_{\mathrm{red}} \to \kappa[\epsilon]/(\epsilon^2)$, which maps $y_i \otimes 1$ onto an element of $\kappa \cdot \epsilon$. (iii) is obvious from this.

(iv) By (iii), $x \in R[\mathcal{U}_a]$ and hence $R[\mathcal{U}_a] = R[x] = R[\widehat{\mathcal{U}}_a]$.

(v) We must have $n_1 > 0$, and by (iii) we also have $n_2 = 0$. It follows that $n_i = 0$ for all $i$ even, and $n_i \leqslant n_1$ for all $i$ odd.

Recall that the affine ring of $(\mathcal{U}_a)_K$ is $K[x]$, and it is a Hopf algebra with comultiplication $\mu : x \mapsto x \otimes 1 + 1 \otimes x$. To ease the notation, we will write $x \otimes 1 = u$, $1 \otimes x = v$. We now examine the condition that $R[\mathcal{U}_a] \subset K[x]$ is closed under comultiplication, i.e. $\mu(\pi^{n_i} x^i) \subset R[\mathcal{U}_a] \otimes_R R[\mathcal{U}_a]$ for all $i$. Clearly,

$$\mu(\pi^{n_i} x^i) = \sum_{j=0}^{i} \pi^{n_i} \binom{i}{j} u^j v^{i-j} \in R[\mathcal{U}_a] \otimes R[\mathcal{U}_a]$$

if and only if

$$n_j + n_{i-j} \leqslant n_i + \mathrm{ord}_K \binom{i}{j},$$

for all $j$. For $i = 2$ and $j = 1$, we get $2n_1 \leqslant e = \mathrm{ord}_K(2)$. For odd $i$ and $j = 1$, we get $n_1 \leqslant n_i$.

Thus the only possible affine ring for $R[\mathcal{U}_a]$ is $R[\pi^c x, x^2]$ with $2 \leqslant 2c \leqslant e$. The assertion on the length of the quotient of the two Lie algebras is easy to check. ∎

**8.7 The constraint on the $\Delta$-invariant** Let $H \in \Xi(\mathcal{G}_K)$ be such that $H \simeq \mathrm{SO}_{2n+1}$, and let $\mathcal{H}$ be the schematic closure of $H$ in $\mathcal{G}$. For simplicity, we now assume that $\mathcal{H} = \mathcal{G}$ so that we can use the notation set up in 8.5.

By the big-cell decomposition for $\mathcal{G}$, $\mathrm{Lie}\,\mathcal{G}_\kappa$ is the direct sum of $\mathrm{Lie}\,\mathcal{S}_\kappa$ and $\mathrm{Lie}(\mathcal{U}_a \otimes \kappa)$ for $a \in \Phi(G, S)$. It follows that $\mathrm{Lie}(\mathcal{G}_\kappa)_{\mathrm{red}}$ is the direct sum of $\mathrm{Lie}\,\mathcal{S}_\kappa$ and $\mathrm{Lie}(\mathcal{U}_a \otimes \kappa)_{\mathrm{red}}$. By Proposition 8.6

(iii) and Corollary 2.4, the weight of $\mathrm{Lie}(\mathcal{U}_a \otimes \kappa)_{\mathrm{red}}$ is either $a$ or $2a$, and hence Proposition 8.6 (iv) or (v) can be applied. Notice that we can combine (iv) and (v) to say that $R[\mathcal{U}_a]$ is of the form $R[\pi^c x, x^2]$ for an integer $c$ such that $0 \leqslant 2c \leqslant \mathrm{ord}_K(2)$.

If $\mathcal{G} = \mathcal{H}$ is smooth over $R$, it is clear that $\Delta_\mathcal{G}(H) = 0$. Now assume that $\mathcal{G}$ is not smooth over $R$. According to 3.3 and Proposition 4.3, $\widehat{\mathcal{G}}_\kappa \to (\mathcal{G}_\kappa)_{\mathrm{red}}$ is a unipotent isogeny. Proposition 8.6 (iv) and Lemma 2.2 show that $\mathcal{U}_a$ is smooth if $a$ is a long root in $\Phi(G, S)$. On the other hand, since the normalizer of $\mathcal{S}(R)$ in $\mathcal{G}(R)$ permutes the short roots of $\Phi(G, S)$ transitively, the above integer $c$ is the same for all the $2n$ short roots $a$.

Again by the big-cell decomposition for $\mathcal{G}$, $\mathrm{Lie}\,\mathcal{G}$ is the direct sum of $\mathrm{Lie}\,\mathcal{S}$ and $\mathrm{Lie}\,\mathcal{U}_a$ for all $a \in \Phi(G, S)$. We also have a similar decomposition of $\mathrm{Lie}\,\widehat{\mathcal{G}}$. This gives us immediately

$$\mathrm{length}_R(\mathrm{Lie}\,\mathcal{G}/\mathrm{Lie}\,\widehat{\mathcal{G}}) = 2nc.$$

Thus the invariant $\Delta_\mathcal{G}(H)$ is simply the integer $c$, with $0 \leqslant 2c \leqslant \mathrm{ord}_K(2)$.

## 9 The existence and uniqueness

**9.1 Some quadratic lattices** We now describe some quadratic lattices relevant to (good) quasi-reductive models of $\mathrm{SO}_{2n+1}$, $n \geqslant 1$. See Lemma 9.2, Theorems 9.3 and 10.4 for their applications.

Let $V = K^{2n+1}$ with standard basis $\{e_{-n}, \ldots, e_{-1}, e_0, e_1, \ldots, e_n\}$, and let $q$ be the quadratic form on $V$ defined by

$$q\left(\sum_{i=-n}^n x_i \cdot e_i\right) = -x_0^2 + \sum_{i=1}^n x_i x_{-i}.$$

Let $G = \mathrm{SO}(q)$, $L_0 = \sum_{i=-n}^n R \cdot e_i$, and $\mathcal{G}_0 = \mathrm{SO}(L_0, q)$ be the schematic closure of $G$ in $\mathrm{GL}(L_0)$. We call a quadratic lattice over $R$ a *Chevalley lattice* if it is isomorphic to $(L_0, u \cdot q)$ for some $n \geqslant 1$, $u \in R^\times$. It is well-known that if $(L', q')$ is a quadratic lattice over $R$ of odd rank, then $\mathrm{SO}(L', q')$ is a Chevalley scheme if and only if $(L', q')$ is a Chevalley lattice.

Fix a $c \in \mathbb{Z}$ such that $0 \leqslant 2c \leqslant \mathrm{ord}_K(2)$. Put $L_c = R\pi^{-c} \cdot e_0 + L_0$ and let $\mathcal{H} = \mathrm{GL}(L_c)$. We say that a quadratic lattice $(L', q')$ over $R$ is *good* if $(L', q')$ is isomorphic to $(L_c, u\pi^{2c} \cdot q)$ for some $0 \leqslant 2c \leqslant \mathrm{ord}_K(2)$, $u \in R^\times$. We say that $(L', q')$ is *potentially good* if $(L', q') \otimes_R R'$ is a good lattice over $R'$ for some local extension $R \subset R'$ of DVRs.

**9.2 Lemma** *For $0 \leqslant 2c \leqslant \mathrm{ord}_K(2)$, the schematic closure $\mathcal{G}$ of $G = \mathrm{SO}(q)$ in $\mathcal{H} = \mathrm{GL}(L_c)$ is a good quasi-reductive model of $G$ whose $\Delta$-invariant equals $c$.*

PROOF. The case of $c = 0$ is well-known so we assume that $1 \leqslant 2c \leqslant \mathrm{ord}_K(2)$. Therefore, $\kappa$ is of characteristic 2 and $\mathrm{ord}_K(2) \geqslant 2$.

**Claim.** The special fiber of $\mathcal{G}$ is non-reduced, and $(\mathcal{G}_\kappa)_{\mathrm{red}}$ is isomorphic to $\mathrm{Sp}_{2n}/\kappa$. The $\Delta$-invariant of $\mathcal{G}$ is $c$.

PROOF. Recall that $\mathcal{G}_0 = \mathrm{SO}(L_0, q)$ is a Chevalley model of $G$. It is clear that the action of $\mathcal{G}_0(\kappa)$ leaves the subspace of $L_0 \otimes \kappa$ generated by $e_0$ invariant. Therefore, $\mathcal{G}_0(R) \subset \mathrm{GL}(L_c)$. It follows that $\mathcal{G}_0(R) = \mathcal{G}(R)$ and we have a natural morphism $\mathcal{G}_0 \to \mathcal{G}$. By Lemma 4.2, $\mathcal{G}$ is a (good) quasi-reductive group scheme of finite type over $R$.

We will first calculate the $\Delta$-invariant of $\mathcal{G}$ when $n = 1$ for clarity. In this case, $\mathcal{H}(R)$ consists of matrices (relative to the basis $\{e_{-1}, e_0, e_1\}$) of the form

$$\begin{pmatrix} x_1 & \pi^c x_2 & x_3 \\ \pi^{-c} y_1 & y_2 & \pi^{-c} y_3 \\ z_1 & \pi^c z_2 & z_3 \end{pmatrix},$$

with $x_i, y_i, z_i \in R$.

Let $T$ be the standard maximal $K$-split torus of $G$, so that $T(A)$ consists of diagonal matrices of the form

$$\lambda(t) = \begin{pmatrix} t & & \\ & 1 & \\ & & t^{-1} \end{pmatrix}, \qquad t \in A^\times$$

for any commutative $K$-algebra $A$. Let $\mathcal{T}$ be the schematic closure of $T$ in $\mathcal{G}$, or equivalently, in $\mathcal{H}$. It is easy to see that $\mathcal{T}$ is simply the Néron-Raynaud model of $T$. In particular, it is smooth.

Let $U_a$ be the root subgroup of $G$ corresponding to the root $a: \lambda(t) \mapsto t$, so that $U_a(A)$ consists of matrices of the form

$$\begin{pmatrix} 1 & 2u & u^2 \\ 0 & 1 & u \\ 0 & 0 & 1 \end{pmatrix}, \qquad u \in A$$

for any commutative $K$-algebra $A$. Let $\mathcal{U}_a$ be the schematic closure of $U_a$ in $\mathcal{G}$ or $\mathcal{H}$. Then the affine ring $R[\mathcal{U}_a]$ of $\mathcal{U}_a$ is a subring of $K[U_a] = K[x]$, and it is generated by $\pi^c x, x^2$, and $(2/\pi^c)x$. As $c \leq \frac{1}{2} \mathrm{ord}_K(2)$, $2/\pi^c \in \pi^c R$, and hence,

$$R[\mathcal{U}_a] = R[\pi^c x, x^2] \subset K[x].$$

Let $v = \pi^c x$, $w = x^2$. Then $R[\mathcal{U}_a] = R[v, w]/(v^2 - \pi^{2c} w)$. Thus we see that the special fiber of $\mathcal{U}_a$ is non-reduced, isomorphic to $\mathbb{G}_a \times \alpha_2$. It is easy to see that $((\mathcal{U}_a)_\kappa)_{\mathrm{red}}$ is in the root subgroup of $(\mathcal{G}_\kappa)_{\mathrm{red}}$ relative to $\mathcal{T}_\kappa$, for the root $2a$.

By 8.7, the $\Delta$-invariant of $\mathcal{G}$ can be calculated by looking at $\mathcal{U}_a$; it is equal to $c > 0$. Therefore, $\mathcal{G}$ is not a Chevalley scheme. So the reduced special fiber $(\mathcal{G}_\kappa)_{\mathrm{red}}$, being the homomorphic image of $(\mathcal{G}_0)_\kappa \simeq \mathrm{SO}_3$ under a non-trivial unipotent isogeny, must be isomorphic to $\mathrm{SL}_2 = \mathrm{Sp}_2$.

In the general case, let $W$ be the subspace spanned by $e_{-1}, e_0, e_1$, and $G' = \mathrm{SO}(q|_W)$. Then the schematic closure $\mathcal{G}'$ of $G'$ in $\mathcal{G}$ (resp. $\mathcal{G}_0$) is the same as that in $\mathrm{GL}(L_c)$ (resp. $\mathrm{GL}(L_0)$), or in $\mathrm{GL}(W \cap L_c)$ (resp. $\mathrm{GL}(W \cap L_0)$), and is what we have studied in the preceding paragraphs. Since the $\Delta$-invariant can be calculated by looking at an $\mathcal{U}_a$ contained in $\mathcal{G}'$ (cf. 8.7), from the $\mathrm{SO}_3$-calculation, we see that $\Delta$-invariant of $\mathcal{G}$ is $c > 0$. It follows that $\mathcal{G}_\kappa$ is non-reduced, $(\mathcal{G}_0)_\kappa \to (\mathcal{G}_\kappa)_{\mathrm{red}}$ is a non-trivial unipotent isogeny, and $(\mathcal{G}_\kappa)_{\mathrm{red}} \simeq \mathrm{Sp}_{2n}$. ∎

**9.3 Theorem** *Let $G$ be a connected reductive algebraic group defined and split over $K$. Let $\Xi$ be the set of normal algebraic subgroups of $G$ which are isomorphic to $\mathrm{SO}_{2n+1}$ for $n \geqslant 1$. Let $\Delta : \Xi \to \mathbb{Z}$ be a function such that $0 \leqslant 2\Delta(H) \leqslant \mathrm{ord}_K(2)$ for all $H \in \Xi$. Then there is a good quasi-reductive model $\mathcal{G}$ of $G$ over $R$ such that $\Delta_{\mathcal{G}}$ is the given $\Delta$ and $\widehat{\mathcal{G}} \simeq \widetilde{\mathcal{G}}$.*

PROOF. Let $\Theta$ be the set of all connected normal almost simple algebraic subgroups of the derived group of $G$. For $H \in \Theta$, let $Z_H$ be the schematic center of $H$. For $H \in \Xi$, let $\mathcal{H}_H$ be a good model of $H$ whose $\Delta$-invariant is $\Delta(H)$. Such a model exists by the preceding lemma.

For $H \in \Theta \smallsetminus \Xi$, let $\mathcal{H}_H$ be a Chevalley model of $H$. Finally, let $S$ be the connected center of $G$ and $\mathcal{S}$ be the Néron-Raynaud model of $S$. Let $F$ be the kernel of the isogeny

$$G' := S \times \prod_{H \in \Theta} H \to G.$$

Let $M$ be any one of $F$, $Z_H$ for $H \in \Theta$. Then $M$ is a group of multiplicative type over $K$ corresponding to a Galois module $\mathrm{Hom}_{\bar{K}}(M, \mathbb{G}_m)$ which is unramified. Such a group has a canonical extension to a $R$-group scheme of multiplicative type.

The canonical extension $\mathcal{Z}_H$ of $Z_H$ is a natural subgroup scheme of $\mathcal{H}_H$ and the canonical extension $\mathcal{F}$ of $F$ is naturally embedded in

$$\mathcal{S} \times \prod_{H \in \Theta} \mathcal{Z}_H.$$

Thus $\mathcal{F}$ is a closed subgroup scheme of

$$\mathcal{G}' := \mathcal{S} \times \prod_{H \in \Theta} \mathcal{H}_H.$$

Let $\mathcal{G} = \mathcal{G}'/\mathcal{F}$ (in the sense of quotients of fppf sheaves [R, Théorème 1 (iv)]). It is easy to see that $\mathcal{G}$ is a good model of $G$ and it is a quasi-reductive $R$-group scheme with $\Delta_{\mathcal{G}} = \Delta$. Moreover, if $\widehat{\mathcal{G}'}$ is the Chevalley scheme which is the smoothening of $\mathcal{G}'$, then $\mathcal{F}$ is also a closed subgroup scheme of $\widehat{\mathcal{G}'}$. Then quotient $\widehat{\mathcal{G}'}/\mathcal{F}$ is also a Chevalley scheme, and is clearly finite over $\mathcal{G}$. Therefore, it is both the smoothening and normalization of $\mathcal{G}$, and we have $\widehat{\mathcal{G}} \simeq \widetilde{\mathcal{G}}$. ∎

**9.4 Theorem** *Let $\mathcal{G}$ and $\mathcal{G}'$ be good quasi-reductive $R$-group schemes of finite type. Then $\mathcal{G}$ is isomorphic to $\mathcal{G}'$ if and only if there is an isomorphism from $\mathcal{G}_K$ to $\mathcal{G}'_K$ inducing a bijection $\xi : \Xi(\mathcal{G}_K) \to \Xi(\mathcal{G}'_K)$ such that $\Delta_{\mathcal{G}} = \Delta_{\mathcal{G}'} \circ \xi$.*

PROOF. The "only if" part is clear. It remains to prove the other implication.

The given condition implies that $\widehat{\mathcal{G}}$ and $\widehat{\mathcal{G}'}$ are Chevalley schemes of the same type, hence there is an isomorphism $f : \widehat{\mathcal{G}} \to \widehat{\mathcal{G}'}$. Let $\lambda$ be the corresponding isomorphism from the function field of $\widehat{\mathcal{G}'}$ to that of $\widehat{\mathcal{G}}$. Let $\mathcal{S}$ be a good torus of $\widehat{\mathcal{G}}$ and $\mathcal{S}' = f(\mathcal{S})$ the corresponding good torus of $\widehat{\mathcal{G}'}$.

We will identify $\Phi = \Phi(\mathcal{G}_K, \mathcal{S}_K)$ with $\Phi(\mathcal{G}'_K, \mathcal{S}'_K)$ via $f$. For $a \in \Phi$, we denote by $U_a$ (resp. $U'_a$) the corresponding root subgroup of $\mathcal{G}_K$ (resp. $\mathcal{G}'_K$), and by $\mathcal{U}_a$ (resp. $\mathcal{U}'_a$) the schematic closure of $U_a$ (resp. $U'_a$) in $\mathcal{G}$ (resp. $\mathcal{G}'$). If $a$ is not a short root in $\Phi(H, \mathcal{S}_K \cap H)$ for some $H \in \Xi(\mathcal{G}_K)$, then we have $\mathcal{U}_a = \widehat{\mathcal{U}}_a \simeq \widehat{\mathcal{U}}'_a = \mathcal{U}'_a$, where the isomorphism is induced by $f$.

Suppose that $a$ is a short root of $\Phi(H, \mathcal{S}_K \cap H)$ for an $H \in \Xi(\mathcal{G}_K)$. Then the affine ring of $\mathcal{U}_a$ can be constructed from (i) the affine ring of $\widehat{\mathcal{U}}_a$, (ii) the action of $\mathcal{S}$ on the affine ring of $\widehat{\mathcal{U}}_a$, and (iii) the function $\Delta$. The same is true for $\mathcal{U}'_a$. Since $\Delta(H) = \Delta(f(H))$, $f$ again induces an isomorphism $\mathcal{U}_a \simeq \mathcal{U}'_a$.

By the big-cell decomposition theorem, $\prod_{a \in \Phi^+} \mathcal{U}_a \times \mathcal{S} \times \prod_{a \in \Phi^-} \mathcal{U}_a$ is canonically isomorphic to an open subscheme $\Omega$ of $\mathcal{G}$. There is a similar open subscheme $\Omega'$ of $\mathcal{G}'$, and there is an isomorphism $\Omega \simeq \Omega'$, whose inducing isomorphism between the function fields of $\Omega$ and $\Omega'$ is simply $\lambda$. Since $\mathcal{G}(R) \to \mathcal{G}(\kappa)$ is surjective, the collection $\{\mathcal{G}_K\} \cup \{g\Omega g^{-1} : g \in \mathcal{G}(R)\}$ is an open cover of $\mathcal{G}$, and we have isomorphisms $\mathcal{G}_K \xrightarrow{\sim} \mathcal{G}'_K$, $g\Omega g^{-1} \xrightarrow{\sim} f(g)\Omega' f(g)^{-1}$, all compatible with $\lambda$. These isomorphisms patch together to give an isomorphism $\mathcal{G} \to \mathcal{G}'$. Thus the theorem is proved. ∎

We observe that this completes the proof of Theorem 8.3, which combines Theorems 9.3 and 9.4.

**9.5 Theorem** *Let $G$ be a connected reductive group defined and split over $K$ and $\widehat{\mathcal{G}}$ the Chevalley model over $R$ of $G$. For any $\Delta : \Xi(G) \to \mathbb{Z}$ such that $0 \leq 2\Delta(H) \leq \operatorname{ord}_K(2)$ for all $H \in \Xi(G)$, there is a unique model $\mathcal{G}_\Delta$ of $G$ which is a good quasi-reductive group scheme of finite type over $R$ with $\mathcal{G}_\Delta(R) = \widehat{\mathcal{G}}(R)$ and $\Delta_{\mathcal{G}_\Delta} = \Delta$. The smoothening of $\mathcal{G}_\Delta$ is isomorphic to the normalization of $\mathcal{G}_\Delta$, which in turn is just the Chevalley model $\widehat{\mathcal{G}}$. Given $\Delta' : \Xi(G) \to \mathbb{Z}$, satisfying the same condition as $\Delta$, there is a morphism $\mathcal{G}_{\Delta'} \to \mathcal{G}_\Delta$ extending the identity morphism on the generic fibers if and only if $\Delta'(H) \leq \Delta(H)$ for all $H \in \Xi(G)$.*

PROOF. The existence part is clear from Theorem 9.3. For the uniqueness, we argue as in the preceding theorem: the big cells of $\mathcal{G}_\Delta$ can be constructed uniquely from the big cells of $\widehat{\mathcal{G}}$ and the function $\Delta$. The uniqueness theorem also shows that the normalization of $\mathcal{G}_\Delta$ is smooth since this is true for the good model provided by Theorem 9.3.

For the last statement, we first notice that from the construction in Theorem 9.3, there is a morphism $\mathcal{G}_{\Delta'} \to \mathcal{G}_\Delta$ if $\Delta'(H) \leq \Delta(H)$ for all $H \in \Xi(G)$. Conversely, if there is a morphism $\mathcal{G}_{\Delta'} \to \mathcal{G}_\Delta$, as in 8.5, we can find a torus $S$ in $G$ whose schematic closures in $\mathcal{G}_\Delta$ and $\mathcal{G}_{\Delta'}$ are good. Form $\Phi = \Phi(G, S)$ and $\{U_a\}_{a \in \Phi}$ with respect to this torus. Let $\mathcal{U}_a$ (resp. $\mathcal{U}'_a$) be the schematic closure of $U_a$ in $\mathcal{G}_\Delta$ (resp. in $\mathcal{G}_{\Delta'}$). Then $\mathcal{G}_{\Delta'} \to \mathcal{G}_\Delta$ induces a morphism $\mathcal{U}'_a \to \mathcal{U}_a$. This together with 8.7 show that $\Delta'(H) \leq \Delta(H)$ for all $H \in \Xi(G)$. ∎

**9.6 Corollary** *Let $\mathcal{G}$ and $\mathcal{G}'$ be quasi-reductive models of $G$ such that there is a morphism $\mathcal{G}' \to \mathcal{G}$ extending the identity morphism on the generic fiber. Then,*

(i) *$\mathcal{G}' \to \mathcal{G}$ is a finite morphism, and hence $\mathcal{G}'(R) = \mathcal{G}(R)$;*

(ii) *$\mathcal{G}$ is good if and only if $\mathcal{G}'$ is good;*

(iii) *Fix $\mathcal{G}$, there are only finitely many quasi-reductive models $\mathcal{G}'$ with a morphism $\mathcal{G}' \to \mathcal{G}$ as above.*

PROOF. (i) It suffices to check this after a base change $R \subset R'$. Therefore, we may assume that both $\mathcal{G}$ and $\mathcal{G}'$ are good (notice that if $\mathcal{G}$ is good, then so is $\mathcal{G} \otimes_R R'$ for any totally ramified local extension $R \subset R'$ of DVRs). Now as $\mathcal{G}(R) = \mathcal{G}(R')$ is hyperspecial, the smoothening $\widehat{\mathcal{G}}$ of $\mathcal{G}$ is also the smoothening of $\mathcal{G}'$. Since the composition $\widehat{\mathcal{G}} \to \mathcal{G} \to \mathcal{G}'$, being the smoothening morphism of $\mathcal{G}'$, is a finite morphism, so is the morphism $\mathcal{G} \to \mathcal{G}'$.

(ii) is clear from (i) and characterization (3) in Lemma 8.1.

(iii) Again, it suffices to verify this after a faithfully flat base change $R \subset R'$. Then we may and do assume that $\mathcal{G}$ is good. The statement then follows from (ii) and the preceding theorem. ∎

**9.7 Completion of the proof of Theorem 1.2 (a)** We proceed as in 5.3 to choose a subset $S$ of a basis $B$ of $A := R[\mathcal{G}]$, and for $I \supset S$, consider $A_I$ and $\mathcal{G}^I$, etc. Again, $\mathcal{G}^I$ is a quasi-reductive model of $G := \mathcal{G}_K$ for each $I \supset S$. Now fix a finite subset $I$ of $B$ containing $S$. By Corollary 9.6, there are only finitely many $R$-subalgebras of $A \otimes_R K$ which contain $A_I$ and correspond to models of $G$ that are quasi-reductive group schemes of finite type over $R$. For any finite subset $J$ of $B$ such that $A_J \supset A_I$, the algebra $A_J$ is one of these. Therefore, the union of $A_J$ for all such $J$, which is simply $A$, is actually a union of finitely many $A_J$'s, and $A$ is hence of finite type over $R$.

## 10 Quasi-reductive models of $SO_{2n+1}$

**The Lie algebra** We retain the notations and hypothesis from 9.1. In particular, the quadratic form $q$ is as in there. We recall that Lie $G$ can be identified with a Lie subalgebra of $\mathrm{End}(V)$, and Lie $\mathcal{G}$ and Lie $\mathcal{G}_0$ are lattices in Lie $G$.

Let $B(-,-)$ be the bilinear form associated to $q$, that is, $B(v,w) = q(v+w) - q(v) - q(w)$. Let $\wedge^2_K V$ be the exterior square of $V$. There is a natural map

$$\iota : \wedge^2_K V \to \mathrm{End}(V), \qquad \iota(a \wedge b) : v \mapsto B(a,v)b - B(b,v)a.$$

**10.1 Lemma** *The map $\iota$ is a $K$-vector space isomorphism from $\wedge^2_K V$ onto Lie $G$. Moreover, the image of $\wedge^2_R L_0$ is Lie $\mathcal{G}_0$, and the image of $\wedge^2_R L_c$ is Lie $\mathcal{G}$.*

PROOF. The first statement is well-known if 2 is invertible in $K$, and can be verified for an arbitrary non-degenerate quadratic form $q$. However, for the particular $q$ we are working with, the statement is true even if $K$ is replaced with $\mathbb{Z}$. In fact, a basis of the Chevalley algebra Lie $G$ is given on pages 192–193 of [B], and one can verify easily that $\iota(e_{-i} \wedge e_i) = H_i$, $\iota(e_i \wedge e_0) = X_{\varepsilon_i}$, $\iota(e_0 \wedge e_{-i}) = X_{-\varepsilon_i}$, $\iota(e_j \wedge e_i) = X_{\varepsilon_i + \varepsilon_j}$, and so on, where $\Phi = \{\pm\varepsilon_i, \pm\varepsilon_i \pm \varepsilon_j\}$ is the root system of $G$ and $\{X_a\}_{a \in \Phi}$ is the Chevalley basis given in [B]. This also shows that the image of $\wedge^2_R L_0$ is Lie $\mathcal{G}_0$.

By Proposition 8.6 and Theorem 8.4, Lie $\mathcal{G}$ is spanned over $R$ by $H_i$, $1 \leqslant i \leqslant n$, $\pi^{-c} X_a$ for short roots $a \in \Phi$, and $X_a$ for long roots $a \in \Phi$. From this, we we see that $\iota$ carries $\wedge^2_R L_c$ onto Lie $\mathcal{G}$. ∎

**10.2 Lemma** *Let $\mathcal{H}$ be a group scheme locally of finite type over a noetherian ring $R$. Let $R \to R'$ be a flat morphism from $R$ to a noetherian ring $R'$. Then the canonical morphism $(\mathrm{Lie}\,\mathcal{H}) \otimes_R R' \to \mathrm{Lie}(\mathcal{H} \otimes_R R')$ is an isomorphism.*

PROOF. By [SGA3, Exp. II, Proposition 3.3 and page 54], Lie $\mathcal{H}$ is $\operatorname{Hom}_R(e^*(\Omega^1_{\mathcal{H}/R}), R)$, where $e$ is the identity section. The lemma follows from this description, the compatibility of the formation of $\Omega^1_{\mathcal{H}/R}$ and base change ([BLR, 2.1, Proposition 3]), and [M2, Theorem 7.11]. ∎

**10.3 Lemma** *Let $V$ be a vector space of dimension $m$ over $K$. Let $d < m$ be a positive integer prime to $m$. Let $R \subset R'$ be a local extension of DVRs and put $K' = \operatorname{Frac} R'$. Let $L'$ be an $R'$-lattice in $V' := V \otimes_K K'$, $M' = \wedge^d_{R'} L'$. Suppose that there is an $R$-lattice $M$ in $\wedge^d_K V$ such that $M \otimes_R R' = M'$. Then there exist an $R$-lattice $L_1$ in $V$, an element $a \in (K')^\times$ such that $L_1 \otimes_R R' = aL'$, and $d \cdot \operatorname{ord}_K(a) \in \mathbb{Z}$, where $\operatorname{ord}_K$ is the valuation on $K'$ normalized so that $\operatorname{ord}_K(K^\times) = \mathbb{Z}$.*

PROOF. We will abbreviate the assumption on $M'$ to "$M'$ descends to a lattice in $\wedge^d_K V$", and so forth. Let $L_0$ be a lattice in $V$ such that $L'_0 := L_0 \otimes_R R'$ is contained in $L'$. Write $L'/L'_0 \simeq \bigoplus_{i=1}^m R'/r_i R'$ and put $M'_0 = \wedge^d_{R'} L'_0$. It follows that $M'/M'_0$ is isomorphic to the direct sum of $R'/(r_{i_1} \cdots r_{i_d})R'$, $i_1 < \cdots < i_d$. By assumption, this implies $\operatorname{ord}_K(r_{i_1} \cdots r_{i_d}) \in \mathbb{Z}$ for all $i_1 < \cdots < i_d$. It follows that we have $\operatorname{ord}_K(r_i) - \operatorname{ord}_K(r_j) \in \mathbb{Z}$ for all $i, j$, and $d \cdot \operatorname{ord}_K(r_i) \in \mathbb{Z}$ for all $i$. We set $a = r_1^{-1}$ and claim that $L'_1 = a.L'$ descends to a lattice in $V$. Notice that $\wedge^d_{R'} L'_1 = a^d M'$ also descends to a lattice in $\wedge^d_K V$. Therefore, we may and do assume that $L' = L'_1$. The only consequence of $L' = L'_1$ that will concern us is that $D' := \wedge^m_{R'}(L')$ descends to a lattice $D$ in $\wedge^m_K V$.

We apply Grothendieck's theory of flat descent [BLR, Chapter 6]. Put $R'' = R' \otimes_R R'$, $K'' = K' \otimes_K K'$. By assumption, there is a canonical descent datum $\varphi_V : V' \otimes_{K',i_2} K'' \to V' \otimes_{K',i_1} K''$, where $i_1, i_2$ are the two natural embeddings of $K'$ into $K''$. There are also canonical descent data $\varphi_M : M' \otimes_{R',i_2} R'' \to M' \otimes_{R',i_1} R''$ and $\varphi_D : D' \otimes_{R',i_2} R'' \to D' \otimes_{R',i_1} R''$. These are compatible in the sense that $\varphi_V$ and $\varphi_M$ induce the same isomorphism $\wedge^d_{K'} V' \otimes_{K',i_2} K'' \to \wedge^d_{K'} V' \otimes_{K',i_1} K''$, and a similarly condition holds for $\varphi_V$ and $\varphi_D$. Our second claim is that $\varphi_V$ restricts to an isomorphism $\varphi_! : L' \otimes_{R',i_2} R'' \to L' \otimes_{R',i_1} R''$. It then follows that $\varphi_!$ satisfies the cocycle condition, and is a descent datum defining a lattice $L = L_1$ which proves the first claim.

By using an $R'$-basis of $L'$, we can regard the datum $\varphi_V$ as an element $g$ of $\operatorname{GL}_m(K'')$, and the compatibility condition with $\varphi_M$ is that its image in $\operatorname{GL}_{m'}(K'')$ lies in $\operatorname{GL}_{m'}(R'')$, where $m' = \binom{m}{d}$ and the morphism $\operatorname{GL}_m \to \operatorname{GL}_{m'}$ is the $d$-th exterior power representation. Similarly the compatibility with $\varphi_D$ shows that image of $g$ under $\det : \operatorname{GL}_m \to \mathbb{G}_m$ lies in $\mathbb{G}_m(R'')$. The second claim is then that $g$ lies in $\operatorname{GL}_m(R'')$, which is now obvious since $\operatorname{GL}_m \to \operatorname{GL}_{m'} \times \mathbb{G}_m$ is a closed immersion of group schemes over $R''$ (or even over $\mathbb{Z}$) since $d$ is prime to $m$. ∎

**10.4 Theorem** *Let $\mathcal{G}$ be a quasi-reductive model of $G = \operatorname{SO}(q)$. Then there exist a unique $\alpha \in \{0, 1\}$ and a unique $R$-lattice $L$ in $V$ such that*

(i) *$\mathcal{G}$ is the schematic closure of $G$ in $\operatorname{GL}(L)$;*

(ii) *$\iota : \wedge^2_K V \simeq \operatorname{Lie} G$ induces an isomorphism $\wedge^2_R L \simeq \pi^\alpha \operatorname{Lie} \mathcal{G}$.*

*Let*

$$c = \operatorname{length}_R(L/L_0) - \frac{\alpha}{2} \cdot \dim V,$$

*where $L_0$ is any Chevalley lattice in V. Then $(L, \pi^{2c-\alpha} \cdot q)$ is potentially good. Moreover, $\mathcal{G}$ is good if and only if $\alpha = 0$ and $(L, \pi^{2c} \cdot q)$ is good, in which case the $\Delta$-invariant of $\mathcal{G}$ is c, and so $2c \leqslant \text{ord}_K(2)$. Conversely, if $(L', q')$ is a good (resp. potentially good) quadratic lattice, then the schematic closure of $\text{SO}(L', q') \otimes_R K$ in $\text{GL}(L')$ is a good quasi-reductive model over R (resp. a quasi-reductive model over R).*

PROOF. First assume that $\mathcal{G}$ is good. In view of the uniqueness assertion of Theorem 9.4, we may and do assume that $\mathcal{G}$ is the model constructed in Lemma 9.2 using the lattice $L_c$. Then by Lemma 10.1, for $L := L_c$, $\iota$ induces an isomorphism $\wedge^2_R L \simeq \text{Lie}\,\mathcal{G}$. We will now prove the uniqueness of L. Indeed, it is easy to show that the only lattices stable under $\mathcal{G}(R)$ are of the form $\pi^a(R\pi^{-b} \cdot e_0 + L_0)$, with $a, b \in \mathbb{Z}$, $0 \leqslant b \leqslant \text{ord}_K(2)$. Among these, only the one with $a = 0, b = c$, i.e. $L = L_c$, satisfies (ii). By definition, $(L, \pi^{2c} \cdot q)$ is good. According to Lemma 9.2, the $\Delta$-invariant of $\mathcal{G}$ equals c.

Now we drop the assumption that $\mathcal{G}$ is good. By Proposition 3.4, there is a local extension $R \subset R'$ of DVRs such that $\mathcal{G}' := \mathcal{G} \otimes_R R'$ is good. Let $K' = \text{Frac}\,R'$. Then by the case we already treated, there is a unique lattice $L'$ in $V' := V \otimes_K K'$ such that $\mathcal{G}'$ is the schematic closure of $G \otimes_K K'$ in $\text{GL}(L')$, and $\iota$ induces an isomorphism $M' = \wedge^2_{R'} L' \simeq \text{Lie}\,\mathcal{G}'$. Let $r \in R'$ be such that $(L', r \cdot q)$ is good.

Let $M \subset \wedge^2_K V$ be the inverse image of $\text{Lie}\,\mathcal{G}$ under $\iota$. Since $\iota$ obviously commutes with base change, by Lemma 10.2, we have $\iota(M \otimes_R R') = \text{Lie}\,\mathcal{G}'$ and hence $M \otimes_R R' = M'$. According to Lemma 10.3, there exist $\alpha \in \{0, 1\}$, $a \in R'$, an R-lattice L in V such that $\text{ord}_K(a) = \alpha/2$ and $L \otimes_R R' \simeq a.L'$.

We claim that this L satisfies (i) and (ii). Indeed, it suffices to check both properties after the base change $R \subset R'$, and these hold by construction. The uniqueness of L is clear. Moreover, $(L, \pi^{2c-\alpha} \cdot q)$ is potentially good since $(a.L', a^{-2}\pi^{2c} \cdot q)$ is good. Finally, the statement $\mathcal{G}$ is good $\iff \alpha = 0$ and $(L, \pi^{2c} \cdot q)$ is good is also clear.

The converse statement follows from Lemma 9.2. ∎

**10.5 Corollary** *Let $\mathcal{G}$ be a quasi-reductive model of $\text{SO}_{2n+1}$ over R. Let $R \subset R'$ be a local extension of DVRs such that $\mathcal{G}' := \mathcal{G} \otimes_R R'$ is good. Then the $\Delta$-invariant of $\mathcal{G}'$ is divisible by e/2, where e is the ramification index of $R'/R$.*

PROOF. We notice that the invariant $c = c(\mathcal{G})$ takes half-integral values. It is designed in such a way that $c(\mathcal{G}') = e \cdot c(\mathcal{G})$ in the context of the corollary. Thus the corollary is obvious. ∎

**10.6 Corollary** *Assume that the c-invariant of $\mathcal{G}$ is an integer and $\mathcal{G}$ admits a good torus, then $\mathcal{G}$ is good.*

PROOF. Let $\mathcal{T}$ be a good torus of $\mathcal{G}$ over R. Then we can define $\mathcal{U}_a$ etc. as in 8.5. Assume that the weight a submodule of $R[\mathcal{U}_a]$ is $Ry$.

Let $R \subset R'$ be a local extension of DVRs such that $\mathcal{G}' := \mathcal{G} \otimes_R R'$ is good. By Proposition 8.6, either $R'[\mathcal{U}_a] = R'[y]$ or $R'[\mathcal{U}_a] = R'[\pi'^{c'} z, z^2]$ for some integer $c' > 0$, where $\pi'$ is a uniformizer of $R'$, and $z \in R'^{\times} \pi'^{-c'} y$. In the former case, $\mathcal{U}_a$ is smooth. From the proof of the preceding corollary, $c'$ is a multiple of the ramification index e of $R'/R$. Let $c = c'/e$. Then $\pi'^{c'}/\pi^c \in R'^{\times}$, and hence, $z \in R'^\times \pi^{-c} y$. Let $x = \pi^{-c} y$. Then $\pi'^{c'} z \in R'^\times \pi^c x$ and $z \in R'^\times x$. So $R'[\mathcal{U}_a] = R'[\pi^c x, x^2]$, which implies that $R[\mathcal{U}_a] = R[\pi^c x, x^2]$. From this it is clear that the normalization of $\mathcal{U}_a$ over $R$ is smooth.

It then follows from the decomposition of the big-cell associated to the good torus $\mathcal{T}$ that the image of $\mathcal{G}(R) \to \mathcal{G}(\kappa)$ contains a big-cell of $\mathcal{G}(\kappa)$, and hence $\mathcal{G}(R) \to \mathcal{G}(\kappa)$ is surjective. Therefore, the conjugates of the big cell by elements of $\mathcal{G}(R)$, together with $\mathcal{G}_K$, form an open covering of $\mathcal{G}$. Again by the decomposition of the big-cell, it is obvious that the smoothening morphism $\widehat{\mathcal{G}} \to \mathcal{G}$ is finite over each member of this covering. Hence, $\widehat{\mathcal{G}} \to \mathcal{G}$ is a finite morphism and so $\widehat{\mathcal{G}}$ is the normalization of $\mathcal{G}$ over $R$. Therefore, $\mathcal{G}$ is good. ∎

**10.7 Examples** (i) As before, let $\pi$ be a uniformizer of $R$ and $b$ an *odd* integer such that $1 \leqslant b \leqslant \operatorname{ord}_K(2)$. Let $q$ be the quadratic form $x^2 + \pi^b yz$ on a rank-3 lattice $L$. Then $(L, q)$ is potentially good. The corresponding model $\mathcal{G}$ has $c$-invariant $b/2 \in \frac{1}{2}\mathbb{Z} \smallsetminus \mathbb{Z}$, therefore, $\mathcal{G}$ is not good. Observe that $\mathcal{G}$ does admit a good torus. It becomes good over a suitable quadratic extension of $R$.

Notice that we can take $R = \mathbb{Z}_2$ (or its maximal unramified extension) and $b = 1$ to get an example of a non-smooth quasi-reductive model of $SO_3$, while the construction in 9.1 doesn't provide any example of a non-smooth quasi-reductive model of $SO_{2n+1}$ over such a DVR.

One may wonder if $c \in \frac{1}{2}\mathbb{Z}/\mathbb{Z}$ is the only obstruction for being good. The following example shows that this is not the case.

(ii) Let $q$ be the quadratic form $-x^2 - \pi y^2 + \pi^2 yz$ in three variables $x, y, z$. Assume that $\pi^2 \mid 2$. Then $q$ is potentially good but not good. Therefore, this gives us a quasi-reductive model $\mathcal{G}$ of $SO_3$ which is not good. The $c$-invariant of $\mathcal{G}$ is 1.

To see this, let $R' = R[\pi']$ with $\pi'^2 = \pi$. Then $q \otimes_R R'$ is $-x'^2 + \pi^2 y' z'$, where $x' = x + \pi' y$, $y' = y$, $z' = z + (2/\pi'^3)x$. Therefore, $q$ is potentially good. On the other hand, if $q$ is $R$-equivalent to $u(-x^2 + \pi^2 yz)$ for some $u \in R^\times$, then $q \bmod \pi^2$ is a multiple of the square of a linear form, which is not the case.

Notice that by Corollary 10.6, $\mathcal{G}$ doesn't have a good torus, even though both $\mathcal{G}_\kappa$ and $\mathcal{G}_K$ contain 1-dimensional tori.

It is an interesting question to classify all potentially good quadratic lattices. By Theorem 10.4, this is equivalent to classifying quasi-reductive models of $SO_{2n+1}$.

# Appendix: Base change and normalization
by Brian Conrad[3], University of Michigan

Let $X$ be a finite-type flat scheme over a Dedekind domain $R$ with fraction field $K$; we write $S$ to denote $\operatorname{Spec}(R)$. Clearly the structure sheaf of $X_{\text{red}}$ is torsion-free over $R$, and so $X_{\text{red}}$ is also flat over $R$. For any finite extension $K'/K$ we shall write $R'$ to denote the integral closure of $R$ in $K'$; this is Dedekind, is semi-local when $R$ is, and induces finite extensions on residue fields.

---

[3]partially supported by an NSF grant and a Sloan fellowship

***A.1* Remark** If $R$ is a henselian DVR then $R'$ is automatically a henselian DVR as well; in particular, $R'$ is again local.

Quite generally, if $R \to R'$ is any extension of Dedekind domains, inducing an extension $K \to K'$ on fraction fields, we write $X'$ to denote $X \otimes_R R'$ and $X'_{\mathrm{red}}$ to denote $(X')_{\mathrm{red}}$ (and not $(X_{\mathrm{red}})'$). The following was proved by Raynaud [An, App. II, Cor. 3], and later by Faltings [dJ, Lemma 2.13], and our aim will be to describe its proof and its relevance to this paper.

***A.2* Theorem (Raynaud–Faltings)** *There exists a finite extension $K'/K$ such that $X'_{\mathrm{red}}$ has geometrically reduced generic fiber and its normalization $\widetilde{X}'$ is $X'$-finite with geometrically normal generic fiber and geometrically reduced special fibers (over $R'$).*

***A.3* Remark** An algebraic scheme $Z$ over a field $k$ is *geometrically reduced* (resp. *geometrically normal*) over $k$ if $Z \otimes_k k'$ is reduced (resp. normal) for any finite inseparable extension $k'/k$, in which case the same is true for any extension field $k'/k$. We will also use the notion of *geometric integrality* over $k$; see [EGA, IV$_2$, §4.5, 6.7.6ff.] for a detailed discussion.

It follows from Serre's homological criteria "$(R_0) + (S_1)$" for reducedness and "$(R_1) + (S_2)$" for normality (of noetherian rings) [M2, pp. 183ff.] that a finite-type flat scheme over a noetherian normal domain is normal if its generic fiber is normal and its other fibers are reduced. Thus, for $K'$ as in Theorem A.2 and any further flat extension $R' \to R''$ to another noetherian normal domain $R''$, the base change $\widetilde{X}' \otimes_{R'} R''$ is normal, and so it is the normalization of the reduced scheme $X'_{\mathrm{red}} \otimes_{R'} R''$. In particular, the normalization of $(X \otimes_R R'')_{\mathrm{red}} = X'_{\mathrm{red}} \otimes_{R'} R''$ is finite over $X \otimes_R R''$ in a uniform sense as we vary $R''/R'$.

The argument of Raynaud uses rigid-geometry and flattening techniques, whereas the argument of Faltings uses the Stable Reduction Theorem for curves. Strictly speaking, Faltings' proof assumes that $R$ is (local and) excellent [M1, Ch. 13], primarily to ensure finiteness of various normalization maps. We shall reduce the general case to the case of complete local $R$ with algebraically closed residue field, and we then use our reduction steps to describe Faltings' method in a way that avoids some technicalities. We first record:

***A.4* Corollary** *Let $G$ be a flat finite-type separated group scheme over a Dedekind domain $R$ with fraction field $K$. There exists a finite extension $K'/K$ such that:*

- *$G'_{\mathrm{red}}$ is a subgroup of $G'$ with smooth generic fiber;*
- *the normalization $\widetilde{G}'_{\mathrm{red}} \to G'_{\mathrm{red}}$ is finite and is a group-smoothening in the sense of [BLR, §7.1];*
- *these properties are satisfied for the extension $R \hookrightarrow R''$ induced by any injective extension of scalars $R' \hookrightarrow R''$ with $R''$ a Dedekind domain.*

PROOF. The geometric generic fiber $G_{\overline{K}}$ has smooth underlying reduced scheme, as it is a group over an algebraically closed field, so by replacing $K$ with a large finite extension we may assume the $S$-flat $G_{\mathrm{red}}$ has smooth generic fiber. It follows that $G_{\mathrm{red}} \times_S G_{\mathrm{red}}$ is reduced, and hence coincides with $(G \times_S G)_{\text{red}}$, so $G_{\text{red}}$ is a subgroup of $G$. Thus, we may rename $G_{\text{red}}$ as $G$ and may assume $G_K$ is smooth, and by Theorem A.2 we may suppose that the normalization $\widetilde{G}$ is $G$-finite with geometrically reduced fibers over $S$, and its formation commutes with Dedekind extension on $R$.

We conclude that $\widetilde{G} \times_S \widetilde{G}$ is $S$-flat with smooth generic fiber and reduced special fibers, so it is normal (by Serre's criterion). The $S$-separatedness of $G$ and the normality (and $S$-flatness) of $\widetilde{G} \times_S \widetilde{G}$ allow us to use the universal property of normalization to construct a group law on $\widetilde{G}$ compatible with the one on its generic fiber $G_K$; note that finiteness of $\widetilde{G}$ over $G$ provides a bijection $\widetilde{G}(S) = G(S)$, so the identity lifts. Since the fibers of the $S$-flat group $\widetilde{G}$ are geometrically reduced, $\widetilde{G}$ is smooth. Thus, by the argument in 3.3, $\widetilde{G}$ is the group-smoothening. ∎

**A.5 Lemma** *There exists a finite extension $K'/K$ and a nonempty open subscheme $U' \subseteq S' = \mathrm{Spec}(R')$ such that $(X_{K'})_{\text{red}}$ is geometrically reduced and the normalization map $\widetilde{X_{U'}} \to (X_{U'})_{\text{red}}$ is finite with connected components of $\widetilde{X_{U'}}$ having geometrically normal and geometrically integral $U'$-fibers.*

This lemma reduces Theorem A.2 to the case of local $R$, since there are only finitely many points in $S' - U'$.

PROOF. Since $(X_{\overline{K}})_{\text{red}}$ is generically smooth, and the nilradical is locally generated by finitely many elements, by chasing $\overline{K}$-coefficients we may find a finite extension $K'/K$ such that $(X'_{\text{red}}) \otimes_{K'} \overline{K}$ is reduced. That is, upon renaming $K'$ as $K$ and renaming $X_{\text{red}}$ as $X$, we may suppose $X_K$ is geometrically reduced. Further coefficient-chasing in allows us to descend the finite normalization $\widetilde{X}_{\overline{K}} \to X_{\overline{K}}$ to a finite (necessarily birational) map $Y_K \to X_K$, at least after extending $K$ a little; since $Y_K$ is geometrically normal over $K$, it is normal and so it is the normalization of $X_K$. Thus, we may assume that the normalization $\widetilde{X}_K$ of $X_K$ is geometrically normal, and moreover (by extending $K$ a little more) that the connected components of $\widetilde{X}_K$ are geometrically integral.

The finite normalization map $\widetilde{X}_K \to X_K$ over the generic fiber of $X$ may be extended to a finite birational map $Y \to X|_U$ with $U \subseteq S$ a dense open, and by shrinking $U$ we may suppose that $Y$ is $U$-flat. Since the connected components of $Y_K = \widetilde{X}_K$ are irreducible, by shrinking $U$ we may suppose that the connected components $Y_i$ of $Y$ are irreducible. Each $Y_{i,K}$ is an irreducible component of $\widetilde{X}_K$, and so is geometrically normal and geometrically integral. By [EGA, IV$_3$, 9.7.8, 9.9.5], there exists a dense open $U_i \subseteq U$ such that each fiber $Y_{i,u}$ is geometrically normal and geometrically integral for all $u \in U_i$. Renaming $\cap U_i$ as $S$, $Y$ has geometrically normal fibers over $S$; thus, $Y$ is normal, and so the finite birational map $Y \to X$ is the normalization. ∎

We may avoid all difficulties presented by the possible failure of normalizations to be finite, via:

**A.6 Theorem** *If $X_K$ is geometrically reduced, then $\widetilde{X} \to X$ is finite.*

PROOF. For any faithfully flat Dedekind extension $R \to R'$ with associated fraction field extension $K \to K'$, $X'$ is $R'$-flat with reduced generic fiber $X_K \otimes_K K'$, so $X'$ is reduced. Thus, $\widetilde{X} \otimes_R R'$ is reduced and is an intermediate cover between $X'$ and its normalization. Since $\widetilde{X}$ is $X$-finite if and only if $\widetilde{X} \otimes_R R'$ is $X'$-finite (as $R \to R'$ is faithfully flat), and the noetherian property of $X'$ ensures that finiteness of its normalization forces finiteness for all intermediate covers, we conclude that it suffices to prove the finiteness of normalization after base change to $R'$. Thus, Lemma A.5 allows us to reduce to the semi-local case, and then we may certainly reduce to the local case. We may then suppose the base is complete, and hence Japanese, so [EGA, IV$_2$, 7.6.5] ensures finiteness of normalizations for finite-type reduced $R$-schemes. ∎

By Lemma A.5 and Theorem A.6, we may assume that $S$ is local, $X$ is normal, and $X_K$ is geometrically normal and geometrically integral over $K$. For any extension $R \to R'$ of Dedekind domains, the base change $X'$ is reduced and its generic fiber $X'_{K'}$ is geometrically normal, so the normalization $\widetilde{X}' \to X'$ is finite (by Theorem A.6). Our problem is to find a finite extension $R'$ such that $\widetilde{X}'$ has geometrically reduced special fibers; keep in mind that $R'$ is usually just semi-local, and not local.

**A.7 Remark** It suffices to prove generic smoothness of special fibers of $\widetilde{X}'$. Indeed, Serre's homological criteria for reducedness and normality ensure that the $R'$-flat normal $\widetilde{X}'$ must have geometrically-reduced special fibers when it has generically-smooth special fibers. This fact allows us to ignore a nowhere-dense closed subset in the special fibers.

**A.8 Definition** Let $R \to R_0$ be a faithfully flat local map between local Dedekind domains, with $K_0/K$ the corresponding extension on fraction fields and $\kappa_0/\kappa$ the extension on residue fields. The extension $R \to R_0$ is *pseudo-unramified* if:

- The maximal ideal $\mathfrak{m}_{R_0}$ is generated by $\mathfrak{m}_R$, and $\kappa_0/\kappa$ is separable algebraic.
- For every finite extension $K'_0/K_0$, there exists a finite extension $K'/K$ such that $K'_0/K_0$ is contained in a $K$-compositum $L$ of $K_0$ and $K'$, and the integral closure of $R_0$ in $L$ is a quotient of $R' \otimes_R R_0$

**A.8.1 Example** If $R$ is local and $R^{\mathrm{h}}$ is its henselization, then $R \to R^{\mathrm{h}}$ is pseudo-unramified. Indeed, the algebraicity of the fraction-field extension $K \to K^{\mathrm{h}}$ ensures that any finite extension of $K^{\mathrm{h}}$ is contained in a $K$-compositum of $K^{\mathrm{h}}$ and a finite extension $K'/K$, and the compatibility on integer rings follows from the more precise statement that $R' \otimes_R R^{\mathrm{h}}$ is the henselization of the semi-local integral extension $R'/R$. This behavior of henselization with respect to integral ring extensions is a special case of [EGA, IV$_4$, 18.6.8].

**A.8.2 Example** If $R$ is local and henselian, with perfect residue field when $K$ has positive characteristic, then the map $R \to \widehat{R}$ to the completion is pseudo-unramified. To see this, let $\widehat{K}'/\widehat{K}$ be a finite extension of the fraction field $\widehat{K}$ of $\widehat{R}$. We may reduce to the cases when this extension is either separable or purely inseparable. The separable case may be settled by using Krasner's lemma to construct a finite separable $K'/K$ such that $K' \otimes_K \widehat{K} \simeq \widehat{K}'$, and then $R' \otimes_R \widehat{R}$ is the completion of $R'$ (due to $R$-finiteness of $R'$ when $K'/K$ is separable).

It remains to treat the purely inseparable case in positive characteristic $p$. We have $\widehat{R} \simeq \kappa[\![y]\!]$, so $\widehat{K} \simeq \kappa(\!(y)\!)$ has a unique inseparable $p^n$-extension, namely $\widehat{K}^{1/p^n} = \kappa(\!(y^{1/p^n})\!) = \widehat{K}(y^{1/p^n})$. This has valuation ring $\widehat{R}[T]/(T^{p^n} - y)$, and we may choose $y$ to be any uniformizer of $\widehat{R}$. Using a uniformizer $y \in R$, we may take $K' = K(y^{1/p^n})$ with integral closure $R' = R[T]/(T^{p^n} - y)$.

The preceding examples allow for further reduction steps in the proof of Theorem A.2, due to:

**A.9 Lemma** *If $R \to R_0$ is pseudo-unramified, it suffices to consider $X \otimes_R R_0$ over $R_0$.*

PROOF. If $K'_0/K_0$ is a finite extension as in Theorem A.2 for $X_0 = X \otimes_R R_0$ over $R_0$, then by slightly increasing $K_0$ we may suppose (by pseudo-unramifiedness) that $K'_0$ is a $K$-compositum of $K_0$ and a finite extension $K'/K$ such that $R'_0$ is a quotient of $R' \otimes_R R_0$. Note that all residue fields at maximal ideals of $R'_0$ are separable algebraic over the corresponding residue fields on $R'$. Pseudo-unramifiedness has done its work, so now replace $K$ and $K_0$ with $K'$ and $K'_0$, and $R$ and $R_0$ with compatible localizations of $R'$ and $R'_0$, and replace $X$ with the (finite) normalization of $X \otimes_R R'$. This reduce us to the case when $X$ is normal with geometrically normal and geometrically integral generic fiber, and the $X_0$-finite normalization $\widetilde{X}_0$ of $X_0$ has geometrically reduced special fibers.

Since $\widetilde{X}_0 \to X_0 = X \otimes_R R_0$ is a finite surjection, each generic point of the special fiber $(X_0)_{s_0}$ of $X_0$ is hit by a generic point of the special fiber of $\widetilde{X}_0$. Since the special fiber of $X_0$ is merely the base change of the special fiber $X_s$ by the extension of residue fields, each generic point $\xi_s$ of $X_s$ is hit by a generic point $\xi'_{s_0}$ of $(\widetilde{X}_0)_{s_0}$ under the canonical map $\widetilde{X}_0 \to X$. Since $\widetilde{X}_0$ and $X$ are finite unions of normal integral schemes, the induced map $\mathcal{O}_{X,\xi_s} \to \mathcal{O}_{X_0,\xi'_{s_0}}$ between local rings is a local extension of DVRs, and hence is faithfully flat. Passing to the quotient by the maximal ideal of $R$ also kills the maximal ideal of $R_0$, so we get a faithfully flat map $\mathcal{O}_{X_s,\xi_s} \to \mathcal{O}_{(X_0)_{s_0},\xi'_{s_0}}$; the target of this map is a field that is linearly disjoint over $\kappa$ with respect to any finite inseparable extension $\kappa'$ of $\kappa$, since $(X_0)_{s_0}$ is $\kappa_0$-smooth near all of its generic points and $\kappa_0$ is separable algebraic over $\kappa$. It follows that the local ring $\mathcal{O}_{X_s,\xi_s}$ of $X_s$ at the generic point $\xi_s$ is also a field that is linearly disjoint from all such $\kappa'$ over $\kappa$. This says that $X_s$ is smooth near $\xi_s$, and since $\xi_s$ was an arbitrary choice of generic point on $X_s$ we conclude that $X_s$ is generically smooth. ∎

By Example A.8.1 and Lemma A.9, we may assume the local base $R$ is henselian. Let $R_0/R$ be a local integral extension with $R_0$ henselian and inducing an algebraic closure $\overline{\kappa}/\kappa$ on residue fields. Assuming Theorem A.2 for $X_0 = X \otimes_R R_0$ over $R_0$, let us deduce it over $R$. There is a finite extension $K'_0/K_0$ such that the normalization $\widetilde{X}'_0$ of $X'_0 = X \otimes_R R'_0$ has geometrically reduced special fiber. Since $K'_0/K_0$ is finite and $K_0/K$ is algebraic, $K'_0$ may be expressed as a $K$-compositum of $K_0$ and a finite extension $K'$ of $K$. By renaming $R'$ as $R$, we may assume that the normalization $\widetilde{X}_0$ of $X_0$ has geometrically reduced fibers. By expressing $R_0$ as a directed union of integral closures $R'$ of $R$ in finite extensions $K'/K$, we see via finiteness of $\widetilde{X}_0 \to X_0$ that there exists such an $R'$ and a finite birational map $Y \to X' = X \otimes_R R'$ that descends $\widetilde{X}_0 \to X_0$. Since $Y \otimes_{R'} R_0 \simeq \widetilde{X}_0$ is normal, so $Y$ becomes normal after a faithfully flat base change ($R' \to R_0$), it follows that $Y$ is normal. Thus, $Y$ is the normalization of $X'$. Renaming $R'$ as $R$ allows us to therefore assume that $\widetilde{X} \otimes_R R_0$ is the normalization $\widetilde{X}_0$ of $X_0$. Since we have already noted that $\widetilde{X}_0$ has geometrically reduced special fiber, and the special fiber of $\widetilde{X}_0 = \widetilde{X} \otimes_R R_0$ is $\widetilde{X}_s \otimes_\kappa \kappa_0$, it follows that $\widetilde{X}$ has geometrically reduced special fiber, as desired.

We may now assume the henselian $R$ has algebraically closed residue field, so by Example A.8.2 and Lemma A.9 we may assume $R$ is complete. More generally, to settle any particular case $X \to \mathrm{Spec}(R)$ over a general local $R$, it is enough to consider the situation after passing to connected components of the normalization of the base-change of $X$ by a local extension $R \to R'$, where $R'$ is a suitable complete DVR with algebraically closed residue field. These reduction steps allow us to use Faltings' proof of [dJ, Lemma 2.13] to prove Theorem A.2. Let us now show how his argument is applied.

*A.10* **Proof of Theorem A.2** As we have explained already, to settle any particular case we may (after suitable finite extension on $K$ and normalization) restrict attention to the case when $R$ is complete with algebraically closed residue field and $X$ is normal with geometrically normal and geometrically integral generic fiber. In particular, $R$ is excellent. The $R$-flatness and the irreducibility of the generic fiber ensure that both fibers of $X$ have the same pure dimension, say $d$, and the application of our reduction steps (if $R$ was originally more general or $X$ was not normal) preserves the hypothesis of the generic fiber having a specified pure dimension $d$. Thus, we may induct on $d$, the case $d = 0$ being trivial.

Suppose $d = 1$. Working locally on $X$, we may assume $X$ is affine and hence quasi-projective, and so by normalizing the projective closure after a suitable finite extension on $K$ (that may possibly be inseparable even if $X_K$ is $K$-smooth), we may assume $X$ is proper with $X_K$ geometrically normal and geometrically integral. Thus, the curve $X_K$ is $K$-smooth. By the Stable Reduction Theorem for curves of genus $\geqslant 2$ [DM] (or see [AW] for the case of an algebraically closed residue field), after a further finite separable extension on $K$ there exists a proper regular $R$-curve $C$ with generic fiber $X_K$ and generically smooth special fiber $C_s$; the same holds for genus $\leqslant 1$ by direct arguments.

Since $R$ is excellent, $C$ is an excellent surface. Since $\kappa$ is algebraically closed, resolution of singularities for excellent surfaces [Ar] and the factorization theorem for proper regular $R$-curves [Ch, Thm. 2.1] ensure that for any two proper normal $R$-curves $Y$ and $Y'$ with the same generic fiber, each generic point on the special fiber $Y_s$ has an open neighborhood in $Y_s$ that is isomorphic to either an open subscheme in $Y'_s$ or an open subscheme in $\mathbf{P}^1_\kappa$; here is is crucial that $\kappa$ is algebraically closed. Consequently, $Y_s$ is generically smooth if and only if $Y'_s$ is generically smooth. Applying this with $Y = X$ and $Y' = C$, we conclude that $X_s$ is generically smooth.

For $d > 1$, we may work locally on $X$ near each generic point of $X_s$, and so may assume $X = \mathrm{Spec}(A)$ is affine with $X_s$ irreducible. We may also (as above) suppose $X$ is normal and $R$ is complete with algebraically closed residue field. Pick a lift $t \in A$ of a $\kappa$-transcendental element in the function field of $X_s$; this defines a dominant $S$-map $\pi : X \to \mathbf{A}^1_S$ that must be flat over the generic point $\eta_s$ of $\mathbf{A}^1_s$, since $\mathcal{O}_{\mathbf{A}^1_S, \eta_s}$ is a DVR. Thus, shrinking $X$ around the generic point of $X_s$ allows us to assume $\pi$ is flat. The localization $\pi_{(\eta_s)} : X \times_{\mathbf{A}^1_S} \mathrm{Spec}(\mathcal{O}_{\mathbf{A}^1_S, \eta_s}) \to \mathrm{Spec}(\mathcal{O}_{\mathbf{A}^1_S, \eta_s})$ is flat with integral generic fiber of dimension $d - 1$, so $\pi_{(\eta_s)}$ has pure relative dimension $d - 1$.

Since $\mathcal{O}_{\mathbf{A}^1_S, \eta_s}$ is a DVR, the induction hypothesis applied to $\pi_{(\eta_s)}$ provides a finite extension $L$ of the fraction field of $\mathcal{O}_{\mathbf{A}^1_S, \eta_s}$, say with $N$ denoting the finite normalization of $\mathrm{Spec}(\mathcal{O}_{\mathbf{A}^1_S, \eta_s})$ in $L$, such that the flat map $(X \times_{\mathbf{A}^1_S} N)_{\mathrm{red}} \to N$ over the semi-local Dedekind $N$ has geometrically reduced generic fiber and has normalization $(X \times_{\mathbf{A}^1_S} N)^\sim$ such that $\widetilde{\pi_N} : (X \times_{\mathbf{A}^1_S} N)^\sim \to N$ has geometrically normal generic fiber and geometrically reduced special fibers.

Let $C \to \mathbf{A}_S^1$ be the finite normalization of $\mathbf{A}_S^1$ in $L$. Since $C$ is a flat normal $R$-curve, we may use the case $d = 1$ to make a finite extension on $K$ so that $C_K$ is geometrically normal and $C_s$ is geometrically reduced. Finiteness of $C$ over $\mathbf{A}_S^1$ ensures that any open subscheme in $C$ containing the generic points of $N_s$ (i.e., the fiber of $C$ over $\eta_s$) contains the preimage of an open subscheme in $\mathbf{A}_S^1$ around $\eta_s$. Thus, to replace $C$ with a sufficiently small open subscheme around $N$, it suffices to replace $\mathbf{A}_S^1$ with a small open subscheme $U$ around $\eta_s$ (and then we replace $C$ and $X$ with $C_U$ and $X_U$; recall that we only need to work generically on $X_s$; see Remark A.7). Since $(X \times_{\mathbf{A}_S^1} C)^\sim \to C$ localized at $N$ is the flat map $\pi_N^\sim$ with geometrically reduced fibers, we may find an open subscheme $U$ around $\eta_s$ so that $(X_U \times_U C_U)^\sim \to C_U$ is flat with geometrically reduced fibers. Thus, we get a flat map $(X_U \times_U C_U)_s^\sim \to (C_U)_s \subseteq C_s$ with geometrically reduced fibers. Thus, geometric reducedness of $C_s$ implies the same for $(X_U \times_U C_U)_s^\sim$.

The finite map $C_U \to U$ localizes at $\eta_s$ to become the finite map $N \to \operatorname{Spec}(\mathcal{O}_{\mathbf{A}_S^1, \eta_s})$ that is flat, so by shrinking $U$ some more around $\eta_s$ we may suppose that $C_U \to U$ is flat. Thus, $X_U \times_U C_U$ is flat over the normal integral $X_U$, and so its normalization $(X_U \times_U C_U)^\sim$ (which is finite over $X_U \times_U C_U$, since the base $S$ is Japanese) is a finite union of integral finite type $S$-flat schemes. Since the map $(X_U \times_U C_U)^\sim \to X_U$ is dominant and finite, hence surjective, and both the source and target have integral connected components, this map must be flat over an open $V \subseteq X_U$ containing the generic point $\xi_s$ of $X_s$ (since $\mathcal{O}_{X, \xi_s}$ is a DVR). Thus, $X_s$ has a dense open subscheme admitting a flat cover by the geometrically reduced scheme $(X_U \times_U C_U)_s^\sim$, and so $X_s$ is generically geometrically reduced. Thus, $X_s$ is generically smooth. ∎

[dJ]     A. J. de Jong: *Smoothness, semi-stability, and alterations*, Publ. Math. IHES, **83**, 51–93 (1996).

[DM]     P. Deligne, D. Mumford: *The irreducibility of the space of curves of given genus*, Publ. Math. IHES, **36**, 75–110 (1969).

[EGA]    A. Grothendieck: *Eléments de géométrie algébrique*, Publ. Math. IHES, **4, 8, 11, 17, 20, 24, 28, 32**.

[SGA3]   A. Grothendieck et al.: *SGA 3: Schémas en Groupes I, II, III*, Lecture Notes in Math. **151, 152, 153**, Springer-Verlag, Heidelberg (1970). http://modular.fas.harvard.edu/sga

[M1]     H. Matsumura: *Commutative algebra* (2nd ed.), Benjamin (1980).

[M2]     H. Matsumura: *Commutative ring theory*, Cambridge University Press, Cambridge (1986).

[MV]     I. Mirković, K. Vilonen: *Geometric Langlands duality and representations of algebraic groups over commutative rings*, preprint (2004).

[R]      M. Raynaud: *Passage au quotient par une relation d'équivalence plate*, Proceedings of a Conference on Local Fields, 78–85, Springer-Verlag (1967).

[V]      A. Vasiu: *Integral canonical models of Shimura varieties of preabelian type*, Asian J. Math. **3**, 401-517 (1999).

[Wa]     W.C. Waterhouse: *Introduction to affine group schemes*, Grad. Text in Math. **66**, Springer-Verlag, New York-Berlin (1979).

[Yu]     J.-K. Yu: *Smooth models associated to concave functions in Bruhat-Tits theory*, preprint (2003).

---

Gopal Prasad                          Jiu-Kang Yu
Department of Mathematics             Department of Mathematics
University of Michigan                Purdue University
Ann Arbor, Michigan 48109             West Lafayette, IN 47907
U.S.A.                                U.S.A.

*Email:* gprasad@umich.edu            *Email:* jyu@math.purdue.edu